\documentclass[journal]{IEEEtran}
\pdfoutput=1
\usepackage{subfigure}
\usepackage{multirow}
\usepackage{hyperref}

\usepackage{color}

\usepackage{cite}
\usepackage[pdftex]{graphicx}

\usepackage[cmex10]{amsmath}
\usepackage{amsfonts}
\usepackage{amssymb}

\usepackage{algorithmic}

\usepackage{float}
\floatstyle{ruled}
\newfloat{algorithm}{htbp}{loa}
\floatname{algorithm}{Algorithm}

\usepackage{multirow}
\usepackage{bm}

\newcommand{\ec}{\em}
\newcommand{\ov}{\overline}
\newcommand{\ul}{\underline}
\newcommand{\sdp}{\succcurlyeq}
\newcommand{\la}{\lambda}
\begin{document}
\title{Distributed Algorithms for Optimal Power Flow Problem}

\author{Albert Y.S. Lam, Baosen Zhang, and David Tse
        
\thanks{The authors are with the Department
of Electrical Engineering and Computer Sciences, University of California, Berkeley,
CA, 94720 USA (e-mail: \{ayslam, zhangbao, dtse\}@eecs.berkeley.edu).}% <-this % stops a space
}

\maketitle

\begin{abstract}
%\boldmath
Optimal power flow (OPF) is an important problem for power generation and it is in general non-convex. With the employment of renewable energy, it will be desirable if OPF can be solved very efficiently so its solution can be used in real time. With some special network structure, e.g. trees, the problem has been shown to have a zero duality gap and the convex dual problem yields the optimal solution. In this paper, we propose a primal and a dual algorithm to coordinate the smaller subproblems decomposed from the convexified  OPF. We can arrange the subproblems to be solved sequentially and cumulatively in a central node or solved in parallel in distributed nodes. We test the algorithms on IEEE radial distribution test feeders, some random tree-structured networks, and the IEEE transmission system benchmarks. Simulation results show that the computation time can be improved dramatically with our algorithms over the centralized approach of solving the problem without decomposition, especially in tree-structured problems. The computation time grows linearly with the problem size with the cumulative approach while the distributed one can have size-independent computation time.
\end{abstract}

% \fi
\IEEEpeerreviewmaketitle

\section{Introduction}
\IEEEPARstart{A}{n electric} power system is the main facility to distribute electricity in modern societies. It is a network connecting power supplies (e.g., thermoelectric generators and turbine steam engines) to consumers. A power grid is generally composed of several subsystems: generation, transmission, substation, distribution, and consumers. The generators generate power which is delivered to the substations at high voltages through the transmission network. The power voltage is stepped down and then distributed to the consumers via the distribution networks. In a typical power system, a few hundreds of generators interconnect to several hundreds of substations. The substations distribute the power to millions of consumers with relatively simpler radial networks with tree-like structures.

In the past, research on power systems mainly focused on the core of the network, i.e., from the generation, via transmission, to the substations. All of the control, planning and optimization was done by a single entity (e.g. an ISO). With the integration of renewal energy and energy storage, self-healing ability, and demand response, the focus is shifted toward the consumer side, i.e. distribution networks, and this new paradigm is called the smart grid \cite{smartgrid}.

The optimal power flow (OPF) is one of the most important problems in power engineering and it aims to minimize the generation cost subject to demand constraints and the network physical constraints, e.g. bus voltage limits, bus power limits, thermal line constraint, etc. Due to the quadratic relations between voltage and power, OPF is non-convex. In general, heuristic approaches have been employed to solve the OPF but they are not guaranteed to yield the optimal solution. To simplify the calculation, with assumptions on lossless power line, constant voltage and small voltage angles, OPF can be linearized and this approximation is also called DC-OPF, which is not accurate under all circumstances \cite{dcopf}. For the complex OPF, \cite{zerogap} suggested solving the problem in its dual form and studied the conditions of the power network with zero duality gap. In \cite{geometry}, it was shown that the duality gap is always zero for network structures such as trees which model distribution networks well. \cite{network_zero}, as an independent work of this paper, decomposes the OPF in terms of cycles and branches and formulates the problem as an second-order cone program for tree networks which is equivalent to that given in \cite{coneprogram}.
In  traditional power systems, OPF is mainly for planning purpose. For example, it is used to determine the system state in the day-ahead market with the given system information. In the smart grid paradigm, due to highly intermittent nature of the renewable, the later the prediction is made, the more reliable it is. If OPF can be solved very efficiently, we may solve the OPF in real time thus mitigating some of the unpredictability.

We aim at solving OPF efficiently. When the system size (e.g. the number of buses) increases, solving the problem in a centralized manner is not practical (this will be verified in the simulation). One possible way is to tackle the problem distributedly by coordinating several entities in the system, each of which handle part of the problem and their collaborative effort solves the whole problem.  To do this, a communication protocol is needed to define what information should be conveyed among the entities. We can learn from the networking protocol development to design a communication protocol for OPF. The earliest form of protocols for the Internet was proposed in 70's. They were designed to handle the increasing volume of traffic sent over the Internet in an ad hoc manner. In 1998, Kelly et al. studied Transmission Control Protocol (TCP), which is one of the core protocols in TCP/IP \cite{kelly}. They showed that TCP can be analyzed with a fundamental optimization problem for rate control and the algorithms developed from the optimization  fit the ad hoc designed variants of TCP. This lays down a framework to design communication protocols for complex systems with reasoning. In this framework, we start with an optimization problem representing the system. By optimization decomposition \cite{nonlinear},  the problem is decomposed into (simpler) subproblems which can be solved by different entities in the system independently. The coordination between the subproblems define the communication protocols (i.e. what and how the data exchange between the entities). \cite{chiang} shows that many problems in communications and networking can be cast under this framework and protocols can be designed through primal and dual decomposition. In this paper, we study OPF by decomposing it into subproblems with primal and dual decomposition. Then we propose the primal and dual algorithms, respectively, to solve OPF in a distributed manner and the algorithms determine the communication protocols. Our algorithms do not assume the existence of a communication overlay with topology different from the power network. In other words, a bus only needs to communicate with its one-hop neighbors in the power network. The algorithms can employed to any power network as long as the strong duality holds. We test the algorithms on IEEE radial distribution test feeders, some random tree-structured networks, and the IEEE transmission system benchmarks. Simulation results show that the algorithms can solve OPF distributively and they are scalable, especially when applied to the distribution network. If we apply the algorithm distributedly, the computational time is independent of the number of buses. If we apply the algorithms (with the problem decomposed) in a central node, the computational time grows linearly with the number of buses.

The rest of this paper is organized as follows. In Section \ref{background}, we give the OPF formulation and the necessary background. Section \ref{algorithms} describes the primal and dual algorithms and the mechanism to recover the optimal voltage from the results of the algorithms. We illustrate the algorithms with two examples in Section \ref{example} and present the simulation results in Section \ref{simulation}. In Section \ref{discussion}, we discuss the characteristics of the algorithms and conclude the paper in Section \ref{conclusion}.

%-----------------------------------------------------------------
\section{Preliminaries} \label{background}
\subsection{Problem Formulation}
Assume that there are $n$ buses in the power network. For buses $i$ and $k$, $i\sim k$ means that they are connected by a power line and $i\not\sim k$ otherwise. Let $z_{ik}$ and $y_{ik}$ be the complex impedance and admittance between $i$ and $k$, respectively, and we have $y_{ik}=\frac{1}{z_{ik}}$. We denote $\textbf{Y}=(Y_{ik},1\leq i,k \leq n)$ as the admittance matrix, where
\begin{align*}
	Y_{ik}=\left\{
	\begin{array}{ll}
		\sum_{l\sim i}{y_{il}}& \text{if } i=k\\
		-y_{ik} & \text{if } i \sim k\\
		0 & \text{if } i \not\sim k.
	\end{array} \right.
\end{align*}
Let $\textbf{v}=(V_1,V_2,\ldots,V_n)^T \in \mathbb{C}^n$ and $\textbf{i}=(I_1,I_2,\ldots,I_n)^T \in \mathbb{C}^n$ be the voltage  and current vectors, respectively. By Ohm's Law and Kirchoff's Current Law, we have $\textbf{i} = \textbf{Y}\textbf{v}$. The apparent power injected at bus $i$ is $S_i= P_i + j Q_i = V_i I_i^H$, where $P_i$ and $Q_i$ are the real and reactive power, respectively, and $H$ means Hermitian transpose. We have the real power vector $\textbf{p}=(P_1,P_2,\ldots,P_n)^T=\text{Re}\{\text{diag}(\textbf{v}\textbf{v}^H \textbf{Y}^H)\}$, where $\text{diag}(\textbf{v}\textbf{v}^H \textbf{Y}^H)$ forms a diagonal matrix whose diagonal is $\textbf{v}\textbf{v}^H \textbf{Y}^H$. We define the cost function of Bus $i$ as ${\text{cost}}_i(P_i)=c_{i2}P_i^2+c_{i1}P_{i}+c_{i0}$, where $c_{i0},c_{i1},c_{i2}\in \mathbb{R}$ and $c_{i2}\geq0, \forall i$. OPF can be stated as
\begin{subequations}
\begin{align}
\text{minimize } 	& \sum_{i=1}^n{{\text{cost}}_i(P_i)}   \label{pdfobj} \\
\text{subject to } & \nonumber \\
 & \underline{V_i} \leq |V_i| \leq \overline{V_i}, \forall i	\label{voltagecons}    \\
 & \underline{P_i}\leq P_i \leq \overline{P}_i, \forall i 			\label{powercons}	\\
 & P_{ik} \leq \overline{P}_{ik}, \forall i,k					\label{linecons}		\\
 & \textbf{p}=\text{Re}\{\text{diag}(\textbf{v}\textbf{v}^H \textbf{Y}^H)\}		\label{physcons}
\end{align}
\end{subequations}				
where $\underline{V_i}$, $\overline{V_i}$, $\underline{P_i}$, $\overline{P_i}$, and $\overline{P}_{ik}$ are the lower and upper voltage limits of bus $i$, the lower and upper power limits of bus $i$, and the real power flow limit between buses $i$ and $k$, respectively.  Eq. (\ref{voltagecons}) is the nodal voltage constraint limiting the magnitude of bus voltage. Eq. (\ref{powercons}) is the nodal power constraint limiting the real power generated or consumed and (\ref{linecons}) is the flow constraint. Eq. (\ref{physcons}) describes the physical properties of the network. In this formulation, $\textbf{p}$ and $\textbf{v}$ are the variables. Eqs. (\ref{powercons}) and (\ref{linecons}) are box constraints with respect to $\textbf{p}$ which are the variables of the objective function (\ref{pdfobj}) and they are relatively easy to handle. Eq. (\ref{voltagecons}) together with (\ref{physcons}) make the problem non-convex and hard to solve.  To illustrate the algorithms, we first consider a simplified version of OPF with $c_{i2}=c_{i0}=0,\forall i$ and neglect (\ref{powercons}) and (\ref{linecons}). Having $c_{i0}=0$ will not affect the optimal solution of the original problem. We will explain how to handle non-zero $c_{i2}$ later. By introducing a $n\times n$ complex matrix $W=(W_{ik},1\leq i,k\leq n)=\textbf{v}\textbf{v}^H$, we can write the simplified OPF in the sequel:
\begin{subequations}
\label{originalprob}
\begin{align}
\text{minimize } 	& \sum_{i=1}^n{c_{i1}P_i}  			 \label{pdfobj2}\\
\text{subject to } & \nonumber \\
 & \underline{V_i}^2 \leq W_{ii} \leq \overline{V_i}^2, \forall i	\label{voltagecons2}    \\
 & \text{rank}(W)=1										\label{rankcons} \\
 & \textbf{p}=\text{Re}\{\text{diag}(\textbf{v}\textbf{v}^H \textbf{Y}^H)\}		\label{physcons2}
\end{align}
\end{subequations}	
Let $\textbf{C}=\text{diag}(c_{11},c_{21},\ldots,c_{n1})$ and $\textbf{M}=(M_{ik},1\leq i,k\leq n)=\frac{1}{2}(\textbf{Y}^H \textbf{C}+\textbf{CY})$. By relaxing the rank constraint (\ref{rankcons}), we have the following semidefinite program (SDP):
\begin{subequations}
\label{optprob}
\begin{align}
\text{minimize } 	& \text{Tr}(\textbf{MW})    \label{pdfobj3}\\
\text{subject to } & \nonumber \\
 & \underline{V_i}^2 \leq W_{ii} \leq \overline{V_i}^2, \forall i	\label{voltagecons3}    \\
 & \textbf{W}\succeq 0		\label{sdpcons}
\end{align}
\end{subequations}

where $\text{Tr}(\cdot)$ is the trace operator. We can solve this SDP at a central control center. However, current algorithms for SDP, e.g. primal-dual interior-point methods \cite{pdipm}, can only handle problems with size up to {\ec several hundreds}. We will decompose the problem into smaller ones by exploring the network structure.

\subsection{Zero Duality Gap} \label{sec:gap}
By \cite{geometry}, the simplified OPF and SDP share the equivalent optimal solution provided that the network has a tree structure, is a lossless cycle, or a combination of tree and cycle. For these kinds of network structures which are typically found in distribution networks, the optimal solution computed from (\ref{optprob}) is exactly the same as that from (\ref{originalprob}). Targeting distribution networks, we can merely focus on (\ref{optprob}). For completeness, the approach in \cite{geometry} is outlined below. 

The dual of \eqref{optprob} is given by 
\begin{align} \label{optdual}
\mbox{maximize } & \sum_{i=1}^n (-\ov{\la}_i \ov{V}_i^2 + \ul{\la}_i \ul{V}_i^2) \\
\mbox{subject to } & \ov{\la}_i \geq 0, \ul{\la}_i \geq 0 \; \forall i \nonumber \\
& \Lambda + M \sdp 0, \nonumber
\end{align}
where $\ov{\la}_i$ and $\ul{\la}_i$ are the Lagrangian multipliers associated with the constraints $W_{ii} \leq \ov{V}_i^2$ and $\ul{V}_i^2 \leq W_{ii}$ respectively. From the KKT conditions, \cite{geometry} showed that \eqref{optprob} always has a solution that is rank 1.     

\subsection{Graph Structure}
We will use the following graph structures to decompose SDP.

Consider a graph $G=(V,E)$, where $V=\{i|1\leq i \leq n\}$ are vertices and $E= \{(i,k)\in V\times V\}$ are edges. Vertices $i$ and $k$ are adjacent if $(i,k)\in E$. A clique $C$ is a subset of $V$ whose induced subgraph is fully connected, i.e., $(i,k)\in E, \forall i,k\in C$. A clique is maximal if it cannot be extended to form a larger one by including any adjacent vertex to the clique. In other words, there does not exist a clique whose proper subset is a maximal clique. A chord is an edge which connects two non-adjacent vertices in a cycle. A graph is chordal if each of its cycles with four or more vertices contains a chord. Thus a chordal graph does not a cycle with four or more vertices. If $G$ is not chordal, we can produce a corresponding chordal graph $\tilde{G}=(V,\tilde{E})$, where $\tilde{E}=E\cup E_f$ and $E_f = \{(i,j)\in V\times V - E\}$ are chords of $G$, called fill-in edges. $\tilde{G}$ is not unique. From $\tilde{G}$, we can compute the set of all possible maximal cliques $\mathcal{C}=\{C_1,\ldots,C_{|\mathcal{C}|}\}$, where $C_i=\{j\in V\}$ whose induced subgraph is complete and maximal. If $G$ is a tree, each pair of vertices connecting by an edge forms a maximal clique. For a tree with $n$ vertices, it can be decomposed into $n-1$ maximal cliques.

For $\textbf{M}$ in (\ref{pdfobj3}), we can induce the corresponding $G$ by having $V=\{i|1\leq i \leq n\}$ and $E=\{(i,k)|M_{i,k}\neq 0\}$. $G$ has a very close relationship with the power network structure because of \textbf{Y}. If all $c_{i1},\forall i$ are non-zero, $G$ directly represents the network.

We use the following procedure to produce $\mathcal{C}$ from $\textbf{M}$.
\begin{enumerate}
\item Construct a graph $G=(V,E)$ from $\textbf{M}$.
\item From $G$, compute Maximum Cardinality Search \cite{MCS} to construct an elimination ordering $\sigma$ of vertices \cite{elimination}.
\item With $\sigma$, perform Fill-In Computation \cite{FIC} to obtain a chordal graph $\tilde{G}$.
\item From $\tilde{G}$, determine the set of maximal cliques $\mathcal{C}$ by the Bron-Kerbosch algorithm \cite{maximalclique}.
\end{enumerate}

Note that similar ideas about maximal cliques have been utilized to develop a parallel IPM for SDP \cite{parallelSDP, sdp1, sdp2}. However, we make use of the ideas to decompose SDP into smaller problems, which can be tackled by any appropriate SDP algorithm, not necessarily IPM. Therefore, our approach is more flexible on that any future efficient SDP algorithm can be incorporated into our framework. 

%------------------------
%------------------------
\section{Algorithms} \label{algorithms}
In this section, we will first present the primal and dual algorithms whose outputs are positive semidefinite matrices. Then we will explain the mechanism to convert such a matrix into the voltage vector.

\subsection{Primal Algorithm} \label{primal}
The objective function (\ref{pdfobj3}) can be expressed as
\begin{align}
	\text{Tr}(\textbf{MW}) = \sum_{i,k=1}^n M_{ik}^H W_{ik}, \label{obj}
\end{align}
where each term $M_{ik}^H W_{ik}$ can be classified into one of the following three categories:
\begin{enumerate}
	\item Ignored terms \\
	Each of which has $M_{ik}=0$. Let $\mathcal{I}=\{(i,k)|M_{ik}=0\}$.
	\item Unique terms \\
	For $M_{ik}\neq 0$, both $i$ and $k$ belong to a unique maximal clique. If $i,k\in C_l$, then $i,k\notin C_r, \forall r\neq l$. Let $\mathcal{U}=\{(i,k)| i,k\in C_l, \forall l, i,k\notin C_r, \forall r\neq l\}$.
	\item Shared terms \\
	For  $M_{ik}\neq 0$, both $i$ and $k$ belongs to more than one maximal clique. %Let $S \subseteq \mathcal{C}$ with $|S|\geq 2$. Suppose it belongs to $S_k=\{C_p\in \mathcal{C}\}$, i.e. $i,j\in C_p,\forall p$.% We denote it $(i,j)\in \phi(S_k)$.	
\end{enumerate}
Then (\ref{obj}) becomes
\begin{align}
	\text{Tr}(\textbf{MW}) =& \sum_{i,k|(i,k)\in \mathcal{I}} M_{ik}^H W_{ik} + \sum_{i,k|(i,k)\in \mathcal{U}-\mathcal{I}} M_{ik}^H W_{ik} \nonumber \\
	&+ \sum_{i,k|(i,k)\notin \mathcal{I} \cup \mathcal{U}} M_{ik}^H W_{ik}, \label{obj1}
\end{align}
where all the ignored terms can be ignored. Since each unique term is unique to each maximal clique, (\ref{obj1}) becomes
\begin{align}
\text{Tr}(\textbf{MW}) =\sum_{\substack{i,k\in C_l,\forall C_l\in \mathcal{C}\\|(i,k)\in \mathcal{U}-\mathcal{I}}} M_{ik}^H W_{ik} + \sum_{i,k|(i,k)\notin \mathcal{I} \cup \mathcal{U}} M_{ik}^H W_{ik}. 
\label{primalobj}
\end{align}

Eq. (\ref{voltagecons3}) gives bounds to each $W_{ii}, 1\leq i\leq n$ and it is equivalent to
\begin{align}
\underline{V}_i^2\leq W_{ii} \leq \overline{V}_i^2, \forall i\in C_l, \forall C_l\in \mathcal{C}.
\end{align}

By \cite{partial}, a matrix is positive semidefinite if all its submatrices corresponding to the maximal cliques induced by the matrix are all positive semidefinite. Let $\textbf{W}_{C_lC_l}$ be the partial matrix of $\textbf{W}$ with rows and columns indexed according to $C_l$. Eq. (\ref{sdpcons}) is equivalent to
\begin{align}
\textbf{W}_{C_lC_l} \succeq 0, \forall C_l\in \mathcal{C}.
\end{align}
Hence, (\ref{optprob}) is written as
\begin{subequations}
\label{optprob2}
\begin{align}
\text{minimize } 	& \sum_{\substack{i,k\in C_l,\forall C_l\in \mathcal{C}\\|(i,k)\in \mathcal{U}-\mathcal{I}}} M_{ik}^H W_{ik} + \sum_{i,k|(i,k)\notin \mathcal{I} \cup \mathcal{U}} M_{ik}^H W_{ik}     \label{pdfobj4}\\
\text{subject to } & \nonumber \\
 &\underline{V}_i^2\leq W_{ii} \leq \overline{V}_i^2, \forall i\in C_l, \forall C_l\in \mathcal{C}	\label{voltagecons4}    \\
 & W_{C_lC_l} \succeq 0, \forall C_l\in \mathcal{C}.		\label{sdpcons4}
\end{align}
\end{subequations}

If we fix all $W_{ik}$ in the shared terms (those in the second summation in (\ref{pdfobj4}), (\ref{optprob2}) can be decomposed into $|\mathcal{C}|$ subproblems, each of which  corresponds to a maximal clique. For $C_l$, we have the subproblem $l$, as follows:
\begin{subequations}
\label{optprob3}
\begin{align}
\text{minimize } 	& \sum_{i,k\in C_l|(i,k)\in \mathcal{U}-\mathcal{I}} M_{ik}^H W_{ik}      \label{pdfobj5}\\
\text{subject to } & \nonumber \\
 &\underline{V}_i^2\leq W_{ii} \leq \overline{V}_i^2, \forall i\in C_l,  i\notin C_r, \forall r\neq l	\label{voltagecons5}    \\
 & \textbf{W}_{C_lC_l} \succeq 0.		\label{sdpcons5}
\end{align}
\end{subequations}

In (\ref{optprob3}), only the semidefinite constraint (\ref{sdpcons5}) involves those variables which are not unique to $C_l$, i.e., $W_{ik}$ such that $(i,k)\notin \mathcal{U}$. By introducing a slack variable $X_{ik,l}=W_{ik}$ for each shared $W_{ik}$ in subproblem $l$, we define $\tilde{\textbf{W}}_{C_lC_l}=(\tilde{W}_{ik},i,k\in C_l)$ and $\tilde{\textbf{M}}_{C_lC_l}=(\tilde{M}_{ik},i,k\in C_l)$ where 
\begin{align*}
\tilde{W}_{ik} = \left\{ \begin{array}{ll}
         {W}_{ik} & \mbox{if $(i,k)\in \mathcal{U}$}\\
        {X}_{ik,l} & \mbox{otherwise}\end{array} \right.
\end{align*}
and
\begin{align*}
\tilde{M}_{ik} = \left\{ \begin{array}{ll}
         {M}_{ik} & \mbox{if $(i,k)\in \mathcal{U}$}\\
        0 & \mbox{otherwise}.\end{array} \right.
\end{align*}
Then (\ref{optprob3}) becomes
\begin{subequations}
\label{optprob4}
\begin{align}
\text{minimize } 	& \text{Tr}(\tilde{\textbf{M}}_{C_lC_l}\tilde{\textbf{W}}_{C_lC_l})      \label{pdfobj6}\\
\text{subject to } & \nonumber \\
 &\underline{V}_i^2\leq W_{ii} \leq \overline{V}_i^2, \forall i\in C_l,  i\notin C_r, \forall r\neq l	\label{voltagecons6}    \\
 & \tilde{\textbf{W}}_{C_lC_l} \succeq 0		\label{sdpcons6}\\
 & X_{ik,l} = W_{ik}, \forall i,k | (i,k)\notin{\mathcal{U}}. \label{slackcons6}
\end{align}
\end{subequations}
Note that $W_{ik}$'s in (\ref{slackcons6}) are given to the subproblem. When given such $W_{ik}$'s, all subproblems are independent and can be solved in parallel. Let the domain of (\ref{optprob4}) be $\Phi_l$. Given $W_{ik}$ where $(i,k)\notin \mathcal{U}$, let $\phi_l(W_{ik}|(i,k)\notin \mathcal{U})=\inf_{\tilde{\textbf{W}}_{C_lC_l}\in \Phi_k}\{\text{Tr}(\tilde{\textbf{M}}_{C_lC_l}\tilde{\textbf{W}}_{C_lC_l})  \}$. (\ref{optprob2}) becomes
\begin{subequations}
\label{optprob5}
\begin{align}
\text{minimize } 	& \sum_{\forall C_l\in \mathcal{C}} \phi_l(W_{ik}|(i,k)\notin \mathcal{U}) + \sum_{i,k|(i,k)\notin \mathcal{U}} M_{ik}^H W_{ik}     \label{pdfobj7}\\
\text{subject to } & \nonumber \\
 &\underline{V}_i^2\leq W_{ii} \leq \overline{V}_i^2, \forall i|(i,i)\notin \mathcal{U}. 	\label{voltagecons7} 
\end{align}
\end{subequations}
Eq. (\ref{optprob5}) is the master problem which minimizes those $W_{ik}$ shared by the maximal cliques. With those shared $W_{ik}$ computed in (\ref{optprob5}), we minimize the $W_{ik}$ unique to each subproblem given in (\ref{optprob4}).

In (\ref{optprob5}), those shared  $W_{ik}$ can be further classified according to nodes and edges:

\begin{enumerate}
	\item Nodes\\
	Let $\lambda_{ii,l}$ be the Lagrangian multiplier for (\ref{slackcons6}) with $i=k$. The subgradient of  $W_{ii}$ with respect to subproblem $l$ is $-\lambda_{ii,l}$ \cite{convex}. Thus the overall subgradient is $\sum_{l|i\in C_l}{(-\lambda_{ii,l})} + M_{ii}^H$. At iteration $t$, we update $W_{ii}$ by
\begin{align}
	W_{ii}^{(t+1)}=Proj\left( W_{ii}^{(t)}-\alpha^{(t)}\left(\sum_{l|i\in C_l}{(-\lambda_{ii,k})} + M_{ii}^H\right) \right), \label{nodecost}
\end{align}
where
\begin{align*}
Proj(x)=\left\{ \begin{array}{ll}
         \underline{V}_i^2 & \mbox{if $x < \underline{V}_i^2$},\\
         \overline{V}_i^2 & \mbox{if $x > \overline{V}_i^2$},\\
        x & \mbox{otherwise},\end{array} \right.
\end{align*}
 $\alpha^{(t)}$ is the step size at iteration $t$, and $W_{ii}^{(t)}$ represents $W_{ii}$ at iteration $t$.
	
	\item Edges in $E$ \\
	We consider (\ref{slackcons6}) with $i\neq k$. Since $W_{ik}$ and $X_{ik,l}$ are complex numbers,  we can handle the real and imaginary parts separately, i.e. $\text{Re}\{X_{ik,l}\}=\text{Re}\{W_{ik}\}$ and $\text{Im}\{X_{ik,l}\}=\text{Im}\{W_{ik}\}$. Let  $\lambda^{\text{Re}}_{ik,l}$ and $\lambda^{\text{Im}}_{ik,l}$ be their corresponding Lagrangian multipliers for subproblem $l$, respectively.  In (\ref{pdfobj6}), $\forall i\neq k$, the $ik$ and $ki$ terms always come in a pair. We have 
\begin{align*}
	M_{ik}^HW_{ik}+M_{ki}^HW_{ki} =& 2\text{Re}\{M_{ik}^H\}\text{Re}\{W_{ik}\} \\
							      & -2\text{Im}\{M_{ik}\}\text{Im}\{W_{ik}\}.
\end{align*}
 A subgradient of the real part of $W_{ik}$ is $\sum_{l|i,k\in C_l}{(-\lambda^{\text{Re}}_{ik,l})} + 2\text{Re}\{M_{ik}^H\}$. At iteration $t$, we update $\text{Re}\{W_{ik}^{(t)}\}$ by
\begin{align}
	\text{Re}\{W_{ik}^{(t+1)}\}=&\text{Re}\{W_{ik}^{(t)}\}- \nonumber\\
					&\alpha^{(t)}\left(\sum_{l|i,k\in C_l}{(-\lambda^{\text{Re}}_{ik,l})}+ 2\text{Re}\{M_{ik}^H\}\right).  \label{edgere}
\end{align}
Similarly, for the imaginary part, we have
\begin{align}
	\text{Im}\{W_{ik}^{(t+1)}\}=&\text{Im}\{W_{ik}^{(t)}\}- \nonumber\\
					&\alpha^{(t)}\left(\sum_{l|i,k\in C_l}{(-\lambda^{\text{Re}}_{ik,l})}- 2\text{Im}\{M_{ik}^H\}\right).  \label{edgeim}
\end{align}
	\item Edges in $E_f$ \\
	Recall that fill-in edges are ``artificial'' edges added to $G$ to make $\tilde{G}$. For $(i,k)\in E_f$, we have $M_{ik}=0,i\neq k$. Similarly, at iteration $k$, we update its real and imaginary parts by
\begin{align}
	\text{Re}\{W_{ik}^{(t+1)}\}=\text{Re}\{W_{ik}^{(t)}\}-\alpha^{(t)}\left(\sum_{l|i,k\in C_l}{(-\lambda^{\text{Re}}_{ik,l})} \right), \label{fedgere}
\end{align} 
and
\begin{align}
	\text{Im}\{W_{ik}^{(t+1)}\}=\text{Im}\{W_{ik}^{(t)}\}-\alpha^{(t)}\left(\sum_{l|i,k\in C_l}{(-\lambda^{\text{Im}}_{ik,l})} \right). \label{fedgeim}
\end{align}	
	
\end{enumerate}

We can interpret the updating mechanism as follows:
certain maximal cliques share a component $W_{ik}$ (if $i=k$, it corresponds to a node; otherwise, it corresponds to an edge or a fill-in edge). $W_{ik}$ represents electricity resources and $-M_{ik}^H$ is its default price. An agent (i.e. a node responsible for computing the update) which is common to all those maximal cliques sharing the resource determines how much resource should be allocated to each maximal clique. In other words, it fixes $W_{ik}$ and every party gets this amount. $X_{ik,l}$ is the actual resource required by $C_l$ and $\lambda_{ik,l}$ corresponds to the price of the resources when $W_{ik}$ is allocated to it. If $C_l$ requires more resource than those allocated, i.e., $X_{ik,l}>W_{ik}$, then $\lambda_{ik,l}>0$.  If the net price, i.e. $\sum_{l|i,k\in C_l}\lambda_{ik,l}-M_{ik}$, is positive, the agent should increase the amount of resource allocating to the maximal cliques because it can earn more. If the net price is negative, then supply is larger than demand and it should reduce the amount of allocated resources.

From (\ref{nodecost})--(\ref{fedgeim}), all shared $W_{ik}$ can be updated independently. The update of each $W_{ik}$ only involves those $\{C_l|C_l\in \mathcal{C}, i,k\in C_l\}$. In other words, (\ref{optprob5}) can be further computed separately according to those maximal cliques shared by each $W_{ik}$.

The pseudocode of the primal algorithm is as follows:
\begin{algorithm}
\caption{Primal Algorithm}
\label{alg:PA}
\begin{tabbing}
 Given $Q, \overline{V}, \underline{V}, \mathcal{C}$\\
1. \=Construct (\ref{optprob4}) for each maximal clique\\
2. \textbf{while} stopping criteria not matched \textbf{do}\\
\>3. \textbf{for} \= each subproblem $l$ (in parallel) \textbf{do}\\
\>\>4. Given $W_{ik}$ with $(i,k)\notin \mathcal{U}$, solve (\ref{optprob4})\\
\>\>5. Return $\lambda_{ik,l} \forall i,k|(i,k)\notin \mathcal{U}$\\
\>6. \textbf{end for}\\
\>7.\= Given $\lambda_{ik,l}\forall l|i,k\in C_l$, update the shared $W_{ik}$ with\\
\>\>(\ref{nodecost})--(\ref{fedgeim}) (in parallel)\\
8. \textbf{end while}
\end{tabbing}
\end{algorithm}

%---------------------------------------------------------------------------------------------------------
\subsection{Dual Algorithm} \label{dual}

Let $\Omega_{ik}=\{C_l|i,k\in C_l, \forall l\}$. Problem (\ref{obj}) can be written as
\begin{align}
&\text{Tr}(\textbf{MW})\nonumber\\
&= \sum_{i,k|(i,k)\in \mathcal{U}} M_{ik}^H W_{ik} + \sum_{i,k|(i,k)\notin \mathcal{U}} |\Omega_{ik}|\frac{M_{ik}^H W_{ik}}{|\Omega_{ik}|}  \nonumber\\
&=\sum_{C_l\in \mathcal{C}}  \left(\sum_{i,k\in C_l|(i,k)\in \mathcal{U}} {M_{ik}^H W_{ik} }+\sum_{i,k\in C_l|(i,k)\notin \mathcal{U}} \frac{M_{ik}^H W_{ik} }{|\Omega_{ik}|}  \right). \label{dualobj1}
\end{align}

Problem (\ref{optprob}) becomes
\begin{subequations}
\label{optprob6}
\begin{align}
\text{minimize } 	&  \sum_{C_l\in \mathcal{C}}  \left(\sum_{\substack{i,k\in C_l\\|(i,k)\in \mathcal{U}}} {M_{ik}^H W_{ik} } +\sum_{\substack{i,k\in C_l\\|(i,k)\notin \mathcal{U}}} \frac{M_{ik}^H W_{ik} }{|\Omega_{ik}|}  \right) 			 \label{pdfobj6}\\
\text{subject to } & \nonumber \\
 & \underline{V}_i^2\leq W_{ii} \leq \overline{V}_i^2, \forall i\in C_l, \forall C_l\in \mathcal{C}	\label{voltagecons6}    \\
 & \textbf{W}_{C_lC_l} \succeq 0, \forall C_l\in \mathcal{C}.			\label{psdcons6} 
\end{align}
\end{subequations}	

Problem (\ref{optprob6}) can be separated into subproblems based on the maximal cliques. However, the subproblems are not completely independent of each other as (\ref{psdcons6}) involves some common variables shared between the subproblems. Similar to the primal algorithm, we can replace $\textbf{W}_{C_lC_l}$ with $\tilde{\textbf{W}}_{C_lC_l}$. For each $W_{ik}|(i,k)\notin \mathcal{U}$, let $X_{ik,l}$ be a copy of $W_{ik}$ in $C_l\in\Omega_{ik}$. To make all $\tilde{\textbf{W}}_{C_lC_l}$ consistent, we should have
\begin{align}
	W_{ik}=X_{ik,l_1}=X_{ik,l_2}=\cdots=X_{ik,l_{|\Omega_{ik}|}}, &\forall l_r|C_{l_r}\in \Omega_{ik}, \nonumber\\
									&\forall W_{ik}|(i,k)\notin \mathcal{U}, \label{slackequality0}
\end{align}
or simply
\begin{align}
	X_{ik,l_1}=X_{ik,l_2}=\cdots=X_{ik,l_{|\Omega_{ik}|}}, &\forall l_r|C_{l_r}\in \Omega_{ik}, \nonumber\\
									&\forall i,k|(i,k)\notin \mathcal{U}. \label{slackequality}
\end{align}
For each $(i,k)\notin \mathcal{U}$, (\ref{slackequality}) can be written into $|\Omega_{ik}|-1$ equalities, e.g.,
%\begin{align*}
%X_{ik,l_1}&=X_{ik,l_2}\\
%X_{ik,l_1}&=X_{ik,l_3}\\
%&\vdots\\
%X_{ik,l_1}&=X_{ik,l_{|\Omega_{ik}|}}
%\end{align*}
%and
\begin{align}
X_{ik,l_1}&=X_{ik,l_2},\nonumber\\
X_{ik,l_2}&=X_{ik,l_3},\nonumber\\
&\vdots \label{equality1}\\
X_{ik,l_{|\Omega_{ik}|-1}}&=X_{ik,l_{|\Omega_{ik}|}}. \nonumber
\end{align}
As shown later, the update mechanism of the dual algorithm depends only on how we arrange (\ref{slackequality}) into equalities. In fact, there are many ways to express the $|\Omega_{ik}|-1$ equalities provided that each $X_{ik,l}$ appears in at least one of the equalities. Suppose the $r$th equality be $\tilde{X}_{ik,r}(1)=\tilde{X}_{ik,r}(2)$. We assign a Lagrangian multiplier  $\upsilon_{ik,r}$ to it. Then we have
\begin{align}
\upsilon_{ik,1}(\tilde{X}_{ik,1}(1)&-\tilde{X}_{ik,1}(2))=0\nonumber\\
\upsilon_{ik,2}(\tilde{X}_{ik,2}(1)&-\tilde{X}_{ik,2}(2))=0\nonumber\\
&\vdots \label{equality2} \\
\upsilon_{ik,|\Omega_{ik}|-1}(\tilde{X}_{ik,|\Omega_{ik}|-1}(1)&-\tilde{X}_{ik,|\Omega_{ik}|-1}(2))=0 \nonumber
\end{align}
When we sum all these equalities up, each $X_{ik,l}$ will be associated with an aggregate Lagrangian multiplier $\tilde{\upsilon}_{ik,l}$, which is composed of all $\upsilon_{ik}$ associated with $X_{ik,l}$. For example, in (\ref{equality1}), we have $\tilde{X}_{ik,1}(1)=X_{ik,l_1}$, $\tilde{X}_{ik,1}(2)=X_{ik,l_2}$, and $\tilde{X}_{ik}(2,1)=X_{ik,l_2}$. Thus $\tilde{\upsilon}_{ik,l_1}=\upsilon_{ik,1}$ and $\tilde{\upsilon}_{ik,l_2}=\upsilon_{ik,2}-\upsilon_{ik,1}$.

In (\ref{equality2}), $\forall 1\leq i,k \leq n$, the corresponding $r$th equality for the $(i,k)$ pair implies the one for the $(k,i)$ pair, i.e.,
\begin{align*}
\tilde{X}_{ik,r}(1)=\tilde{X}_{ik,r}(2) \Rightarrow \tilde{X}_{ki,r}(1)=\tilde{X}_{ki,r}(2),
\end{align*}
due the positive semidefinite property of $\textbf{W}$ given in (\ref{sdpcons}).
We have
\begin{align*}
&\upsilon_{ik,r}\tilde{X}_{ik,r}(1) = \upsilon_{ik,r}\tilde{X}_{ik,r}(2) \\
&\Rightarrow (\upsilon_{ik,r}\tilde{X}_{ik,r}(1))^H = (\upsilon_{ik,r}\tilde{X}_{ik,r}(2))^H\\
 &\Rightarrow \upsilon_{ik,r}^H \tilde{X}_{ik,r}^H(1) = \upsilon_{ik,r}^H \tilde{X}_{ik,r}^H(2)\\
 &\Rightarrow \upsilon_{ik,r}^H \tilde{X}_{ki,r}(1) = \upsilon_{ik,r}^H \tilde{X}_{ki,r}(2).
\end{align*}
Thus the aggregate Lagrangian multiplier for $X_{ki,l}$ can be computed directly from that for $X_{ik,l}$, i.e., $\tilde{\upsilon}_{ki,l}=\tilde{\upsilon}_{ik,l}^H$.

Let $\tilde{\upsilon}=(\tilde{\upsilon}_{ik,l_r},i,k\in C_l|(i,k)\notin\mathcal{U}, \forall C_l\in \mathcal{C}; l_r | C_{l_r}\in\Omega_{ik})$. We form the dual function $d(\tilde{\upsilon},\textbf{W})$ by aggregating (\ref{slackequality}) into (\ref{dualobj1}). We have
\begin{align}
&d(\tilde{\upsilon},\textbf{W}) \nonumber \\ 
&=\sum_{C_l\in \mathcal{C}}  \left(\sum_{i,k\in C_l|(i,k)\in \mathcal{U}} {M_{ik}^H W_{ik} }+\sum_{i,k\in C_l|(i,k)\notin \mathcal{U}} \frac{M_{ik}^H W_{ik} }{|\Omega_{ik}|}  \right) \nonumber\\
&+\sum_{i,k|(i,k)\notin \mathcal{U}}\sum_{r=1|C_{l_r}\in \Omega_{ik}}^{|\Omega_{ik}|}{\tilde{\upsilon}_{ik,l_r}X_{ik,l_r}} \nonumber\\
&=\sum_{C_l\in \mathcal{C}}  \left(\sum_{\substack{i,k\in C_l\\|(i,k)\in \mathcal{U}}} {M_{ik}^H W_{ik} }+\sum_{\substack{i,k\in C_l\\|(i,k)\notin \mathcal{U}}} \left(\frac{M_{ik}^H W_{ik} }{|\Omega_{ik}|} +\tilde{\upsilon}_{ik,l}X_{ik,l}  \right) \right)\nonumber\\
&\triangleq \sum_{C_l\in \mathcal{C}} d(\tilde{\upsilon},\tilde{\textbf{W}}_{C_lC_l})
\end{align}
Given $\tilde{\upsilon}$, (\ref{optprob6}) becomes
\begin{subequations}
\label{optprob7}
\begin{align}
\text{minimize } 	&  \sum_{C_l\in \mathcal{C}} d(\tilde{\upsilon},\tilde{\textbf{W}}_{C_lC_l}) \label{pdfobj7}\\
%\sum_{C_l\in \mathcal{C}}  \left(\sum_{i,k\in C_l|(i,k)\in \mathcal{U}} {M_{ik}^H W_{ik} }\right. \nonumber\\
%				&\left.+\sum_{i,k\in C_l|(i,k)\notin \mathcal{U}} \left(\frac{M_{ik}^H W_{ik} }{|\Omega_{ik}|} +\tilde{\upsilon}_{ik,l}X_{ik,l}  \right) \right) 			 \label{pdfobj7}\\
\text{subject to } & \nonumber \\
 & \underline{V}_i^2\leq W_{ii} \leq \overline{V}_i^2, \forall i\in C_l, \forall C_l\in \mathcal{C}	\label{voltagecons7}    \\
 & \tilde{\textbf{W}}_{C_lC_l} \succeq 0, \forall C_l\in \mathcal{C}.			\label{psdcons7} 
\end{align}
\end{subequations}
Then problem (\ref{optprob7}) can be divided into subproblems according to the maximal cliques and each of them is independent of each other. Subproblem $l$ is stated as:
\begin{subequations}
\label{sdsdp1}
\begin{align}
\text{minimize } 	&   \sum_{\substack{i,k\in C_l\\|(i,k)\in \mathcal{U}}} {M_{ik}^HW_{ik}} +\sum_{\substack{i,k\in C_l\\|(i,k)\notin \mathcal{U}}} \left(\frac{M_{ik}^H}{|\Omega_{ik}|}+\tilde{\upsilon}_{ik,l}  \right)X_{ik,l} \\
\text{subject to } & \nonumber \\
 & \underline{V}_i^2\leq W_{ii} \leq \overline{V}_i^2, \forall i\in C_l|(i,i)\in\mathcal{U}, \\
 & \underline{V}_i^2\leq X_{ii,l} \leq \overline{V}_i^2, \forall i\in C_l|(i,i)\notin\mathcal{U}, \\
 & \tilde{\textbf{W}}_{C_lC_l} \succeq 0.	
\end{align}
\end{subequations}

By solving (\ref{sdsdp1}), we denote the optimal $\tilde{X}_{ik,r}(z)$ for the $r$th equality in (\ref{equality2}) by $\tilde{X}_{ik,r}^{opt}(z)$, where $z\in\{1,2\}$. The gradient of $-d(\tilde{\upsilon},\tilde{\textbf{W}}_{C_lC_l})$\footnote{In the dual form, we maximize $\inf_{\tilde{\textbf{W}}}d(\tilde{\upsilon},\tilde{\textbf{W}})$ over $\tilde{\upsilon}$. In minimization, we consider $-d(\tilde{\upsilon},\tilde{\textbf{W}}_{C_lC_l})$.} with respect to $\upsilon_{ik,r}$ is
\begin{align*}
-\frac{\partial d}{\partial \upsilon_{ik,r}}=\tilde{X}_{ik,r}^{opt}(2)-\tilde{X}_{ik,r}^{opt}(1).
\end{align*}
Let $\upsilon_{ik,r}^{(t)}$, $\alpha^{(t)}>0$, and $\tilde{X}_{ik,l}^{(t)}(z)$  be the Lagrangian multiplier of the $r$th equality associated with $W_{ik}$, the step size, $\tilde{X}_{ik,l}^{opt}(z)$, respectively, at time $t$. Then we can update $\upsilon_{ik,r}$ in (\ref{equality2}) by
\begin{align}
\upsilon_{ik,r}^{(t+1)} = \upsilon_{ik,r}^{(t)} - \alpha^{(t)}\left( \tilde{X}_{ik,r}^{(t)}(2)-\tilde{X}_{ik,r}^{(t)}(1)\right). \label{dualpriceupdate}
\end{align}

If $\tilde{X}_{ik,r}^{(t)}(2)<\tilde{X}_{ik,r}^{(t)}(1)$, then $\upsilon_{ik,r}^{(t+1)}>\upsilon_{ik,r}^{(t)}$. This will make the coefficient corresponding to $\tilde{X}_{ik,r}(1)$ larger  while making that corresponding to $\tilde{X}_{ik,r}(2)$ smaller. At time $t+1$, the subproblem will obtain $\tilde{X}_{ik,r}^{(t+1)}(1)<\tilde{X}_{ik,r}^{(t)}(1)$ and $\tilde{X}_{ik,r}^{(t+1)}(2)>\tilde{X}_{ik,r}^{(t)}(2)$. Hence, $|\tilde{X}_{ik,r}^{(t+1)}(2)-\tilde{X}_{ik,r}^{(t+1)}(1)|<|\tilde{X}_{ik,r}^{(t)}(2)-\tilde{X}_{ik,r}^{(t)}(1)|$. 
On the other hand, if $\tilde{X}_{ik,r}^{(t)}(2)>\tilde{X}_{ik,r}^{(t)}(1)$, then $\upsilon_{ik,r}^{(t+1)}<\upsilon_{ik,r}^{(t)}$. This will make the coefficient corresponding to $\tilde{X}_{ik,r}(1)$ smaller  while making that corresponding to $\tilde{X}_{ik,r}(2)$ larger. Then we will get $\tilde{X}_{ik,r}^{(t+1)}(1)>\tilde{X}_{ik,r}^{(t)}(1)$ and $\tilde{X}_{ik,r}^{(t+1)}(2)<\tilde{X}_{ik,r}^{(t)}(2)$. This will also make $|\tilde{X}_{ik,r}^{(t+1)}(2)-\tilde{X}_{ik,r}^{(t+1)}(1)|<|\tilde{X}_{ik,r}^{(t)}(2)-\tilde{X}_{ik,r}^{(t)}(1)|$. Therefore, (\ref{dualpriceupdate}) drives $X_{ik,l}$'s in (\ref{slackequality}) become closer to each other in value when the algorithm evolves. In other words, (\ref{dualpriceupdate}) tries to make equality (\ref{slackequality}) hold when the algorithm converges.

At any time before the algorithm converges, i.e., (\ref{slackequality}) does not hold, the solution $\textbf{W}$ with the computed $X_{ik,l},\forall i,k| (i,k)\notin\mathcal{U},\forall l|C_l\in \Omega_{ik}$ is an infeasible solution. We can always construct a feasible $\hat{\textbf{W}}$ with $W_{ik}$ which is the average of all $X_{ik,l}$ in (\ref{slackequality}).

The purpose of (\ref{dualpriceupdate}) is to make the two entity $\tilde{X}_{ik,r}^{(t)}(1)$ and $\tilde{X}_{ik,r}^{(t)}(2)$ closer to each other. As long as $\tilde{X}_{ik,r}^{(t)}(1)$ and $\tilde{X}_{ik,r}^{(t)}(2)$ have been computed (from two subproblems), we can update $\upsilon_{ik,r}$ with \ref{dualpriceupdate}. Thus different $\upsilon_{ik}$ can be updated asynchronously.  Since the only co-ordination between subproblems is through (\ref{dualpriceupdate}), synchronization is not required in dual algorithm. 

We can interpret the updating mechanism as follows: $\upsilon_{ik,r}$ is the price assigned to equality $\tilde{X}_{ik,r}(1)=\tilde{X}_{ik,r}(2)$. We can treat $\tilde{X}_{ik,r}(1)$ and $\tilde{X}_{ik,r}(2)$ as demand and supply of electricity resources, respectively. If the demand is larger than the supply, i.e., $\tilde{X}_{ik,r}(1)>\tilde{X}_{ik,r}(2)$, we should increase the price so as to suppress the demand and to equalize the supply and demand. On the other hand, if the supply is larger than the demand, i.e. $\tilde{X}_{ik,r}(1)<\tilde{X}_{ik,r}(2)$, we should reduce the price in order to boost the demand.

The pseudocode of the dual algorithm is as follows:
\begin{algorithm}
\caption{Dual Algorithm}
\label{alg:DA}
\begin{tabbing}
 Given $Q, \overline{V}, \underline{V}, \mathcal{C}$\\
1. \=Pair up slack variables for the shared variables into equalities\\
2. Construct (\ref{sdsdp1}) for each maximal clique\\
3. \textbf{while} stopping criteria not matched \textbf{do}\\
\>4. \textbf{for} \= each subproblem $l$ (in parallel) \textbf{do}\\
\>\>5. Given $\tilde{\upsilon}_{ik}$, solve (\ref{sdsdp1})\\
\>\>6. Return $X_{ik,l}, \forall i,k|(i,k)\notin \mathcal{U}$\\
\>7. \textbf{end for}\\
\>8.\= Given $X_{ik,l_r}$, update the price $\upsilon_{ik,l_r}$ with (\ref{dualpriceupdate}) (in parallel\\
\>\> and asynchronously)\\
9. \textbf{end while}
\end{tabbing}
\end{algorithm}

\subsection{Computation of Voltage}
When either the primal or the dual algorithm converges, assuming zero duality gap, we obtain the optimal $\textbf{W}=\textbf{vv}^H$. To obtain each bus voltage and voltage flown on each line, we first compute the voltage magnitude at each bus, $|V_i|=\sqrt{W_{ii}},1\leq i\leq n$. For $1\leq i,k \leq n$, if there is a line between nodes $i$ and $k$, the corresponding line voltage angle difference $\theta_{ik}$ can be found by solving $W_{ik}=|V_i||V_k|e^{j\theta_{ik}}$ at either bus $i$ or $k$. For the former, bus $k$ needs to send $|V_k|$ to bus $i$, and vice versa. By fixing the voltage angle of a particular bus to zero, the voltages of the whole network can be found subsequently.

\section{Quadratic Cost Function}
Up to now we have focused on the OPF problem with a linear objective function. In practice, sometimes a quadratic cost function is used. If this is the case, the methods developed so far can be used as subroutines to solve the OPF problem by adding a outer loop to the iteration.  

Let $\text{cost}_i (P_i)=c_{i2} P_i^2 + c_{i1} P_i$ be the cost function associated with $P_i$. We assume this function is convex for all buses, that is, $ c_{i2} > 0 \; \forall i$. From \eqref{physcons}, $P_i= \text{Tr}(\mathbf{A_i vv^H})$, where $\mathbf{A_i}=\frac12( \mathbf(Y^H E_i)+\mathbf{E_i Y})$ and $\mathbf{E_i}$ is the matrix with $1$ in the $(i,i)$th entry and zero everywhere else. Now the OPF problem is (compare with \eqref{originalprob}) 
\begin{subequations}
\label{quadratic}
\begin{align}
\text{minimize } 	& \sum_{i=1}^n{c_{i2} \text{Tr}(A_i W)^2 + c_{i1} \text{Tr}(A_i W)}  			 \label{pdf}\\
\text{subject to } & \nonumber \\
 & \underline{V_i}^2 \leq W_{ii} \leq \overline{V_i}^2, \forall i	\label{voltagecons2}    \\
 & \text{rank}(W)=1										\label{rankcons} \\
 & \textbf{p}=\text{Re}\{\text{diag}(\textbf{v}\textbf{v}^H \textbf{Y}^H)\}		\label{physcons2}
\end{align}.
\end{subequations}
Using Schur's complement, we may write \eqref{quadratic} equivalently as 
\begin{subequations}
\label{Schur}
\begin{align}
\text{minimize } & \sum_{i=1}^n (t_i + c_{i1} \text{Tr} (A_i W)  \\
\text{subject to } & \begin{bmatrix} 
t_i & \sqrt{c_{i2}} \text{Tr}(\mathbf{A_i W}) \\
\sqrt{c_{i2}} \text{Tr}(\mathbf{A_i W}) & 1 \end{bmatrix}\succeq 0 \; \forall i \label{schur} \\
& \underline{V_i}^2 \leq W_{ii} \leq \overline{V_i}^2, \forall i \nonumber \\
& \text{rank}(W)=1. \nonumber
\end{align}
\end{subequations}
Relax the the rank $1$ constraint and taking the dual, we get
\begin{subequations}
\label{dualquad}
\begin{align} 
\mbox{maximize } & \sum_{i=1}^n (-\ov{\la}_i \ov{V}_i^2 +\ul{\la}_i \ul{V}_i^2 - u_i) \\
\mbox{subject to } & \sum_{i=1}^n (c_{i1} A_i - 2\sqrt{c_{i2}}z_i A_i)+\Lambda \sdp 0  \\
& \begin{bmatrix}
1 & z_i \\
z_i & u_i \end{bmatrix} \sdp 0 \; \forall i \label{dualz},   
\end{align}
\end{subequations}
where the constraint \eqref{dualz} corresponds to the Schur's compliment constraint in \eqref{Schur}. The constraint \eqref{dualz} can be rewritten as $u_i \geq z_i^2$, for a given $z_i$, the maximizing $u_i$ is $z_i^2$. Therefore the we may drop the constraints \eqref{dualz} and replace the $u_i$ in the objective function by $z_i^2$. If we fix the $z_i$'s, then \eqref{dualquad} becomes a function of $z_i$'s  
\begin{align} \label{dualJz}
J(\mathbf{z})= 
\mbox{maximize } & \sum_{i=1}^n (-\ov{\la}_i \ov{V}_i^2 +\ul{\la}_i \ul{V}_i^2 - z_i^2) \\
\mbox{subject to } & \sum_{i=1}^n (c_{i1} A_i - 2\sqrt{c_{i2}}z_i A_i)+\Lambda \sdp 0 \nonumber. 
\end{align}
For fixed $\mathbf{z}$, \eqref{dualJz} is in the form of \eqref{optdual}. Therefor $J(\mathbf{z})$ is a dual of the optimization problem with linear cost functions with costs $(c_{11}-2\sqrt{c_{12}}z_1,c_{21}-2\sqrt{c_{22}}z_2,\dots,c_{n1}-2\sqrt{c_{n2}}z_n)$. To find the optimal solution of \eqref{dualJz} we may use any of the algorithm in the previous sections. Let $W^*(\mathbf{z})$ denote the optimal solution to $J(\mathbf{z})$. To find the optimal $\mathbf{z}$, we use a gradient algorithm. The Lagrangian of \eqref{dualJz} is 
\begin{align}
\mathcal{L}(\mathbf{\la},W)=& \sum_{i=1}^n (-\ov{\la}_i \ov{V}_i^2 +\ul{\la}_i \ul{V}_i^2 - z_i^2)\nonumber \\
& + \text{Tr}((\sum_{i=1}^n (c_{i1} A_i - 2\sqrt{c_{i2}}z_i A_i)+\Lambda)W).
\end{align}
By a standard result in convex programming, the gradient of $J(z)$ is given by $\frac{J(\mathbf{z})}{z_i}=\frac{\partial \mathcal{L}(\mathbf{\la}^*,W^*)}{\partial z_i}= -2 \text{Tr} ( \sqrt{c_{i2}} A_i W^*)-2 z_i$, where $(\mathbf{\la}^*,W^*)$ is a pair of optimal dual-primal solutions (dependent on $\mathbf{z}$). Therefore, to solve the problem with quadratic cost functions, we add an additional outer loop to the solution algorithms for the linear cost functions.

%---------------------------------------------------------------------------------------- 
\section{Illustrative Examples} \label{example}
In this section, we will consider two examples to illustrate how the algorithms work. The first one is a $n$-bus radial network with height equal to one while the other is a five-bus network with a four-bus ring. The former gives ideas how many maximal cliques share a common bus. The latter demonstrates multiple-level bus and edge sharing. We will also show where the different components are implemented and how the communication is accomplished.

\begin{figure}[!t]
	\begin{center}
		\subfigure[Toplogy]{\label{fig:tree}\includegraphics[width=4cm]{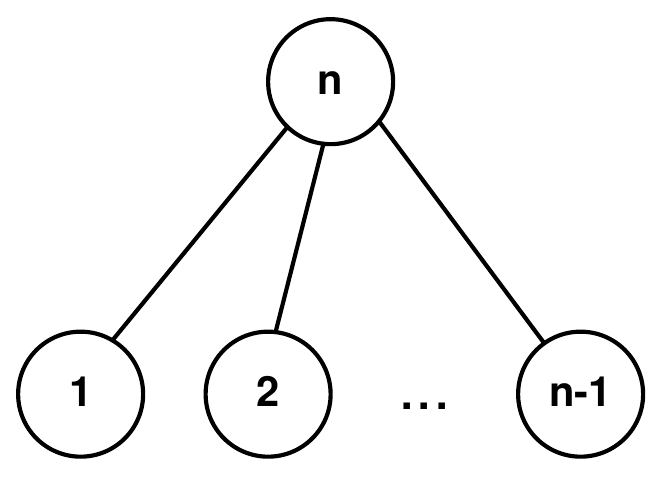}}
    		\subfigure[Communications between nodes $l$ and $n$]{\label{fig:tree_comm}\includegraphics[width=4cm]{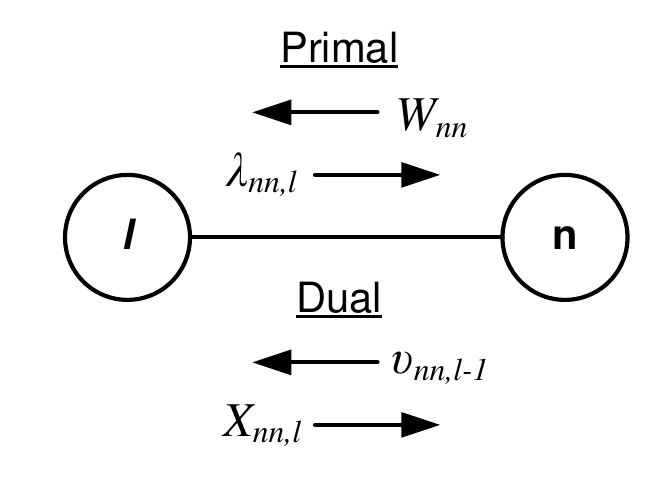}} \\
	\end{center}
	\caption{$n$-bus radial network}
  \label{treenetwork}
\end{figure}

\subsection{$n$-Bus Radial Network}
Consider the topology of the network given in Fig.~\ref{fig:tree}. We have
\begin{align*}
\textbf{M}=
 \begin{pmatrix}
  M_{11} 	& 0		 & 0		&	\cdots & M_{1n} \\
 0		& M_{22} 	 & 0		&	\cdots & M_{2n} \\
 0		& 0 		 & M_{33} &	\cdots & M_{3n} \\
  \vdots  & \vdots	 & \vdots   & \ddots 	  & \vdots  \\
  M_{n1} & M_{n2}      & M_{n3} &  \cdots & M_{nn}
 \end{pmatrix}.
\end{align*}

\subsubsection{Primal Algorithm}
Eq. (\ref{primalobj}) becomes
\begin{align*}
\text{Tr}(\textbf{MW})
&=\underbrace{\sum_{l=1}^{n-1}{(M_{nl}^{H}W_{nl}+M_{ln}^{H}W_{ln}+M_{ll}^{H}W_{ll})}}_{\text{unique terms}} +\underbrace{M_{nn}^{H}W_{nn}.}_{\text{shared term}}
\end{align*}
Each branch with the end nodes forms a maximal clique. We have $\mathcal{C}=\{C_l|1\leq l \leq n-1\}$ where $C_l=\{l,n\}$ and $\mathcal{U}={(n,n)}$. By introducing a slack variable $X_{nn,l}$ for $W_{nn}$ to $C_l$, we have
\begin{align*}
\tilde{\textbf{W}}_{C_lC_l} =
 \begin{pmatrix}
  W_{ll} & W_{ln}  \\
  W_{nl} & X_{nn,l}  
 \end{pmatrix}.
\end{align*}
Let
\begin{align*}
\textbf{M}_l =
 \begin{pmatrix}
  M_{ll} & M_{ln}  \\
  M_{nl} &0  
 \end{pmatrix}.
\end{align*}
Given $W_{nn}$, subproblem $l$ for $C_l$ is stated as
\begin{align}
\label{ex1prisubp}
\begin{array}{rc}
	\text{minimize} 	& \text{Tr}(\textbf{M}_l\tilde{\textbf{W}}_{C_lC_l})    \\ 
	\text{subject to}	& \underline{V}_l^2 \leq \text{Tr}\begin{pmatrix}1 & 0\\ 0 & 0\end{pmatrix}\tilde{\textbf{W}}_{C_lC_l} \leq \overline{V}_l^2	    \\
				&  \text{Tr}\begin{pmatrix}0 & 0\\ 0 & 1\end{pmatrix}\tilde{\textbf{W}}_{C_lC_l} = W_{nn}	    \\
				& \tilde{\textbf{W}}_{C_lC_l} \succeq 0
		\end{array}
\end{align}
which is an SDP with a $2\times 2$ variable. Recall that $\phi_l(W_{nn})$ is the optimal value of subproblem $l$ given $W_{nn}$. The master problem is
\begin{align*}
\begin{array}{rc}
	\text{minimize} 	& \sum_{l=1}^{n-1} \phi_l(W_{nn}) + M_{nn}^{H}W_{nn}    \\ 
	\text{subject to}	& \underline{V}_n^2 \leq W_{nn} \leq \overline{V}_n^2.	
		\end{array}
\end{align*}
$\lambda_{nn,l}$ is the Lagrangian multiplier for equality $X_{nn,l}=W_{nn}$ of subproblem $l$. We update $W_{nn}$ by
\begin{align}
W_{nn}^{(t+1)}=\text{Proj}\left(W_{nn}^{(t)}-\alpha^{(t)}\left(\sum_{l=1}^{n-1}{(-\lambda_{nn,l})} + M_{nn}^H\right)\right)
\label{ex1priupdate}
\end{align}
in  $[\underline{V}_n^2,\overline{V}_n^2]$.

In iteration $t$, bus $n$ broadcasts $W_{nn}^{(t)}$ to bus $l,1\leq l \leq n-1$. After receiving $W_{nn}^{(t)}$, For all $l$, bus $l$ solves its own subproblem (\ref{ex1prisubp}) by any suitable SDP method, e.g. primal-dual IPM \cite{pdipm}, in parallel and then returns $\lambda_{nn,l}$ to bus $n$.\footnote{For most of the interior-point methods, primal and dual solutions come in pair. When the algorithm finds the optimal $\tilde{\textbf{W}}_{C_lC_l}$, it will also give $\lambda_{nn,l}$. Thus no extra calculation is required to determine $\lambda_{nn,l}$.} After receiving all $\lambda_{nn,l}$, node $n$ updates $W_{nn}$ by (\ref{ex1priupdate}). The communication pattern is shown in Fig. \ref{fig:tree_comm}.

\subsubsection{Dual Algorithm}

Eq. (\ref{dualobj1}) becomes
\begin{align*}
\text{Tr}(MW)=\sum_{l=1}^{n-1}{(M_{nl}^{H}W_{nl}+M_{ln}^{H}W_{ln}+M_{ll}^{H}W_{ll}+\frac{M_{nn}^{H}W_{nn}}{n-1})}
\end{align*}
As only $W_{nn}$ is common to all maximal cliques, we have $\Omega_{nn}=\{C_1,\ldots,C_{n-1}\}$ and
\begin{align}
X_{nn,1}=X_{nn,2}=\cdots=X_{nn,n-1}=W_{nn}. \label{ex1equality}
\end{align}
Assume that (\ref{ex1equality}) is arranged as follows. We assign Lagrangian mulipliers to the equalities and we have
\begin{align}
\upsilon_{nn,1}(X_{nn,1}&-X_{nn,2})=0\nonumber\\
\upsilon_{nn,2}(X_{nn,1}&-X_{nn,3})=0\nonumber\\
	      &\vdots \label{ex1dualpairup}\\
\upsilon_{nn,n-2}(X_{nn,1}&-X_{nn,n-1})=0.\nonumber
\end{align}
For $1\leq l \leq n-2$, $\tilde{X}_{nn,l}(1)=X_{nn,1}$ and $\tilde{X}_{nn,l}(2)=X_{nn,l+1}$. Then we have
\begin{align*}
\tilde{\upsilon}_{nn,1}&=\upsilon_{nn,1}+\cdots+\upsilon_{nn,n-2},\\
\tilde{\upsilon}_{nn,l}&=-\upsilon_{nn,l-1},  2\leq l \leq n-1.
\end{align*}
%The dual function becomes
%\begin{align*}
%d(\tilde{\upsilon},\textbf{W}) =& \sum_{l=1}^{n-1}\left[M_{nl}^{H}W_{nl}+M_{ln}^{H}W_{ln}+M_{ll}^{H}W_{ll}\right.\\
%&+\left.\left(\frac{M_{nn}^{H}}{n-1}+\tilde{\upsilon}_{nn}(l-1)\right)X_{nn,l}\right]\\
%\end{align*}
Let
\begin{align*}
\textbf{M}_l =
 \begin{pmatrix}
  M_{ll} & M_{ln}  \\
  M_{nl} & \tilde{M}_l  
 \end{pmatrix}
\end{align*}
where 
\begin{align*}
\tilde{M}_l = \left\{ \begin{array}{ll}
         \left(\frac{M_{nn}}{n-1}+\upsilon_{nn,1}+\cdots+\upsilon_{nn,n-2}\right), & l=1\\
        \left(\frac{M_{nn}}{n-1}-\upsilon_{nn,l-1}\right), & 2\leq l\leq n-1\end{array} \right.
\end{align*}
and
\begin{align*}
\tilde{\textbf{W}}_{C_lC_l} =
 \begin{pmatrix}
  W_{ll} & W_{ln}  \\
  W_{nl} & X_{nn,l}  
 \end{pmatrix}.
\end{align*}
Subproblem $l$ is stated as
\begin{align}
\label{ex1dualsubp}
\begin{array}{rc}
	\text{minimize} 	& \text{Tr}(\textbf{M}_l\tilde{\textbf{W}}_{C_lC_l})    \\ 
	\text{subject to}	& \underline{V}_l^2 \leq  \text{Tr}\begin{pmatrix}1 & 0\\ 0 & 0\end{pmatrix}\tilde{\textbf{W}}_{C_lC_l} \leq \overline{V}_l^2	    \\
	& \underline{V}_n^2 \leq  \text{Tr}\begin{pmatrix}0 & 0\\ 0 & 1\end{pmatrix}\tilde{\textbf{W}}_{C_lC_l} \leq \overline{V}_n^2	    \\
				& \tilde{\textbf{W}}_{C_lC_l} \succeq 0
		\end{array}
\end{align}
which is also an SDP with a $2\times 2$ variable. We update the price by, for $1\leq r \leq n-2$,
\begin{align}
\upsilon_{nn,r}^{(t+1)} = \upsilon_{nn,r}^{(t)} - \alpha^{(t)}\left( X_{nn,r+1}-X_{nn,1}\right). 
\label{ex1dualupdate}
\end{align}
At time $t$, node $l$ from $C_l, 1\leq l \leq n-1$, solves (\ref{ex1dualsubp})\footnote{In fact, any bus in a maximal clique can be elected to solve the subproblem. In this case, node $l$ is chosen to reduce the computation concentrated at bus $n$.}, e.g., by IPM, in parallel and sends $X_{nn_l}$ from the optimal $\textbf{W}_{C_lC_l}$ to bus $n$. Whenever any pair  of $X_{nn}$ specified in (\ref{ex1dualpairup}) (e.g., $X_{nn,1}$ and $X_{nn,l}$ for $C_1$ and $C_{l}$, respectively) reach bus $n$, price $\upsilon_{nn,l-1}$ can be updated with (\ref{ex1dualupdate}) by bus $n$ and the updated $\upsilon_{nn,l-1}$ is multicast back to the corresponding buses (e.g., nodes $1$ and $l$). The communication pattern is shown in Fig. \ref{fig:tree_comm}.

\subsection{Five-Bus Network with a Four-Bus Ring}

\begin{figure}[!t]
	\begin{center}
		\subfigure[Primal algorithm]{\label{fig:ring_primal}\includegraphics[width=5cm]{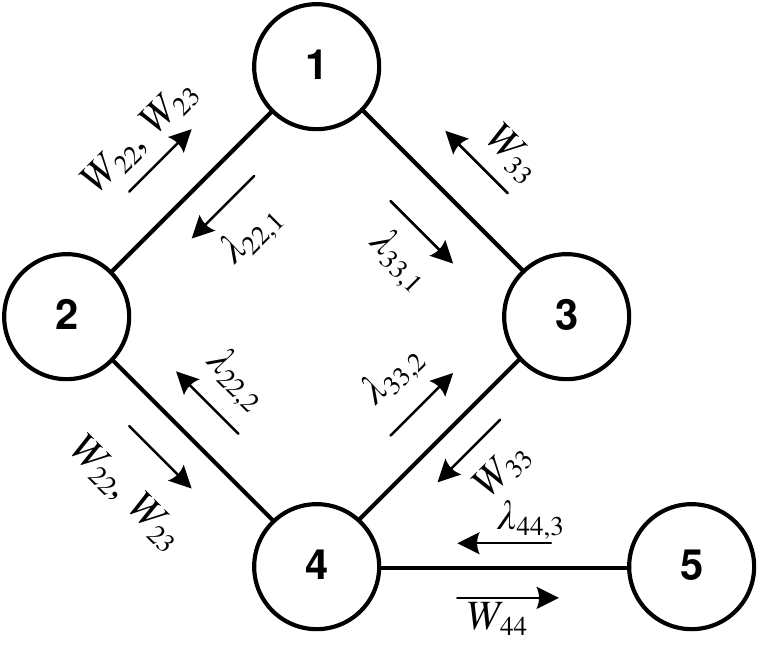}}
    		\subfigure[Dual algorithm]{\label{fig:ring_dual}\includegraphics[width=5cm]{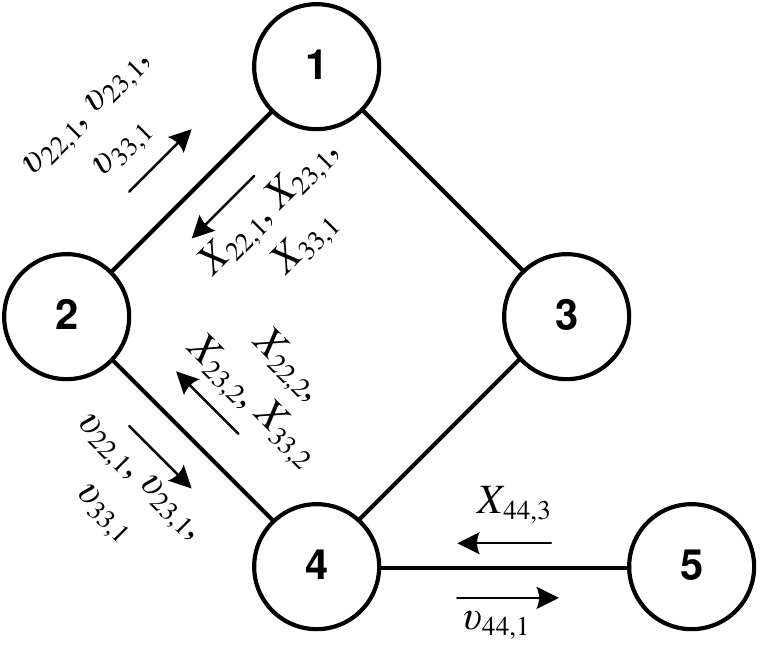}} \\
	\end{center}
	\caption{Structure and communication patterns for the five-bus network with a four-bus ring}
  \label{example2}
\end{figure}

Consider the topology of the network given in Fig.~\ref{fig:ring_primal}. We have
\begin{align*}
\textbf{M}=
 \begin{pmatrix}
  M_{11} 	& M_{12}	& M_{13}	&	0 	& 0 \\
 M_{21}	& M_{22} 	 & 0		& M_{24} 	& 0 \\
 M_{31}	& 0 		 & M_{33} & M_{34} & 0 \\
  0  & M_{42}	 & M_{43}   & M_{44} 	  & M_{45}  \\
  0 & 0      & 0 &  M_{54} & M_{55}
 \end{pmatrix}
\end{align*}

To make it chordal,  suppose we add a fill-in edge between buses 2 and 3 and we have $C_1=\{1,2,3\}$, $C_2=\{2,3,4\}$, and $C_3=\{4,5\}$. Assume bus 2 is used to co-ordinate $C_1$ and $C_2$, and bus 4 for $C_2$ and $C_3$.

\subsubsection{Primal Algorithm}
Let 
\begin{align*}
\footnotesize 
\begin{array}{cc}
\tilde{\textbf{W}}_{C_1C_1} = \begin{pmatrix}
  W_{11} & W_{12} & W_{13}  \\
  W_{21} & X_{22,1} & X_{23,1}\\
  W_{31} & X_{23,1} & X_{33,1}  
 \end{pmatrix},&
 \textbf{M}_1=
  \begin{pmatrix}
  M_{11} & M_{12} & M_{13} \\
  M_{21} & 0 	      & 0 \\
  M_{31} & 0 		&0  
 \end{pmatrix},\\
\tilde{\textbf{W}}_{C_2C_2} =
 \begin{pmatrix}
  X_{22,2} & X_{23,2} & W_{24}  \\
  X_{32,2} & X_{33,2} & W_{34}\\
  W_{42} & W_{43} & X_{44,2}  
 \end{pmatrix},&
  \textbf{M}_2=
 \begin{pmatrix}
  0		& 0		 & M_{24} \\
  0 	 	& 0 	      	& M_{34} \\
  M_{42} 	& M_{43} 	&0  
 \end{pmatrix},\\
\tilde{\textbf{W}}_{C_3C_3} =
 \begin{pmatrix}
  X_{44,3} & W_{45}  \\
  W_{54} & W_{55}   
 \end{pmatrix}, &
  \textbf{M}_3=
 \begin{pmatrix}
  0		& M_{45} \\
  M_{54} 	& M_{55} 
 \end{pmatrix}
\end{array}
\end{align*}

We have
\begin{itemize}
\item Subproblem 1: Given $W_{22}$ and $W_{33}$,
\begin{align}
\label{ex2prisubp1}
\begin{array}{rc}
	\text{minimize} 	& \text{Tr}(\textbf{M}_1\tilde{\textbf{W}}_{C_1C_1})    \\ 
	\text{subject to}	& \underline{V}_1^2 \leq W_{11} \leq \overline{V}_1^2	    \\
				&  X_{ii,1} = W_{ii}, i=2,3	    \\
				& \tilde{\textbf{W}}_{C_1C_1} \succeq 0.
		\end{array}
\end{align}
\item Subproblem 2: Given $W_{22}$, $W_{33}$, and $W_{44}$,
\begin{align}
\label{ex2prisubp2}
\begin{array}{rc}
	\text{minimize} 	& \text{Tr}(\textbf{M}_2\tilde{\textbf{W}}_{C_2C_2})    \\ 
	\text{subject to}	& X_{ii,2} = W_{ii}, i=2,3,4	    \\ 
				& \tilde{\textbf{W}}_{C_2C_2} \succeq 0.
		\end{array}
\end{align}
\item Subproblem 3: Given $W_{44}$,
\begin{align}
\label{ex2prisubp3}
\begin{array}{rc}
	\text{minimize} 	& \text{Tr}(\textbf{M}_3\tilde{\textbf{W}}_{C_3C_3})    \\ 
	\text{subject to}	& \underline{V}_5^2 \leq W_{55} \leq \overline{V}_5^2	    \\ 
				& X_{44,3} = W_{44}	    \\ 
				& \tilde{\textbf{W}}_{C_3C_3} \succeq 0.
		\end{array}
\end{align}
\end{itemize}

In iteration $t$, bus 2 sends $W_{22}$ and $W_{23}$ to both buses 1 and 4 and bus 3 sends $W_{33}$ to both buses 1 and 4.  Bus 4 sends $W_{44}$ to bus 5. Buses 1, 4, and 5 compute (\ref{ex2prisubp1}), (\ref{ex2prisubp2}), and (\ref{ex2prisubp3}), respectively. Then bus 1 sends $\lambda_{22,1}$ to bus 2 and $\lambda_{33,1}$ to bus 3. Bus 4 sends $\lambda_{22,2}$ to bus 2 and $\lambda_{33,2}$ to bus 3. Bus 5 sends $\lambda_{44,3}$ to bus 4, which has $\lambda_{44,2}$.
Buses 2 and 3 update $W_{ii}$ by
\begin{align}
W_{ii}^{(t+1)}=\text{Proj}\left(W_{ii}^{(t)}-\alpha^{(t)}\left(-\lambda_{ii,1}-\lambda_{ii,2} + M_{ii}^H\right)\right)
\label{ex2priupdate1}
\end{align}
within $[\underline{V}_i^2,\overline{V}_i^2]$, for $i=2,3$, respectively. Bus 4 updates $W_{44}$ by
\begin{align}
W_{44}^{(t+1)}=\text{Proj}\left(W_{44}^{(t)}-\alpha^{(t)}\left(-\lambda_{44,2}-\lambda_{44,3} + M_{44}^H\right)\right)
\label{ex2priupdate2}
\end{align}
within $[\underline{V}_4^2,\overline{V}_4^2]$.
Bus 2 updates $W_{23}$ by
\begin{align}
	\text{Re}\{W_{23}^{(t+1)}\}=\text{Re}\{W_{23}^{(t)}\}-\alpha^{(t)}(-\lambda^{\text{Re}}_{23,1}-\lambda^{\text{Re}}_{23,2}), \label{ex2priupdate3}
\end{align} 
and
\begin{align}
	\text{Im}\{W_{23}^{(t+1)}\}=\text{Im}\{W_{23}^{(t)}\}-\alpha^{(t)}(-\lambda^{\text{Im}}_{23,1}-\lambda^{\text{Im}}_{23,2}). \label{ex2priupdate4}
\end{align}
The communication pattern is shown in Fig. \ref{fig:ring_primal}.

\subsubsection{Dual Algorithm}
Assume the equalities for the slack variables are arranged as follows. With Lagrangian multipliers, we have
\begin{align}
\upsilon_{22,1}(X_{22,1}&-X_{22,2})=0\nonumber\\
\upsilon_{23,1}(X_{23,1}&-X_{23,2})=0\nonumber\\
\upsilon_{33,1}(X_{33,1}&-X_{33,2})=0 \label{ex1dualpairup}\\
\upsilon_{44,1}(X_{44,2}&-X_{44,3})=0.\nonumber
\end{align}
Let
\begin{align*}
 \textbf{M}_1=
  \begin{pmatrix}
  M_{11} & M_{12} & M_{13} \\
  M_{21} & \frac{M_{22}}{2}+\upsilon_{22,1}	& \frac{M_{23}}{2}+\upsilon_{23,1}   \\
  M_{31} & \frac{M_{32}}{2}+\upsilon_{23,1}^H 		&\frac{M_{33}}{2}+\upsilon_{33,1} 
 \end{pmatrix},
\end{align*}
\begin{align*}
 \textbf{M}_2=
  \begin{pmatrix}
  \frac{M_{22}}{2}-\upsilon_{22,1}	& \frac{M_{23}}{2}-\upsilon_{23,1}  & M_{24} \\
  \frac{M_{32}}{2}-\upsilon_{23,1}^H &\frac{M_{33}}{2}-\upsilon_{33,1} 	& M_{34}\\
  M_{42}						& M_{43}						& \frac{M_{44}}{2}+\upsilon_{44,1}
 \end{pmatrix},
\end{align*} and
\begin{align*}
 \textbf{M}_3=
  \begin{pmatrix}
  \frac{M_{44}}{2}-\upsilon_{44,1} 	& M_{45}\\
  M_{54}						& M_{55}
 \end{pmatrix}.
\end{align*}
We have
\begin{itemize}
\item Subproblem 1: Given $\upsilon_{22,1}$, $\upsilon_{23,1}$, and $\upsilon_{33,1}$,
\begin{align}
\label{ex2dualsubp1}
\begin{array}{rc}
	\text{minimize} 	& \text{Tr}(\textbf{M}_1\tilde{\textbf{W}}_{C_1C_1})    \\ 
	\text{subject to}	& \underline{V}_i^2 \leq W_{ii} \leq \overline{V}_i^2, i=1,2,3	    \\
				& \tilde{\textbf{W}}_{C_1C_1} \succeq 0.
		\end{array}
\end{align}
\item Subproblem 2: Given  $\upsilon_{22,1}$, $\upsilon_{23,1}$, $\upsilon_{33,1}$, and $\upsilon_{44,1}$,
\begin{align}
\label{ex2dualsubp2}
\begin{array}{rc}
	\text{minimize} 	& \text{Tr}(\textbf{M}_2\tilde{\textbf{W}}_{C_2C_2})    \\ 
	\text{subject to}	& \underline{V}_i^2 \leq W_{ii} \leq \overline{V}_i^2, i=2,3,4	    \\ 
				& \tilde{\textbf{W}}_{C_2C_2} \succeq 0.
		\end{array}
\end{align}
\item Subproblem 3: Given $\upsilon_{44,1}$,
\begin{align}
\label{ex2dualsubp3}
\begin{array}{rc}
	\text{minimize} 	& \text{Tr}(\textbf{M}_3\tilde{\textbf{W}}_{C_3C_3})    \\ 
	\text{subject to}	&  \underline{V}_i^2 \leq W_{ii} \leq \overline{V}_i^2, i=4,5	    \\ 
				& \tilde{\textbf{W}}_{C_3C_3} \succeq 0.
		\end{array}
\end{align}
\end{itemize}

At time $t$, bus 2 announces $\upsilon_{22,1}$, $\upsilon_{23,1}$, and $\upsilon_{33,1}$ to buses 1 and 4,  and bus 4 sends $\upsilon_{44,1}$ to bus 5. Buses 1, 4, and 5 solve (\ref{ex2dualsubp1}), (\ref{ex2dualsubp2}), and (\ref{ex2dualsubp3}), respectively. Then bus 1 sends $X_{22,1}$, $X_{23,1}$, and $X_{33,1}$ to bus 2. Bus 4 sends $X_{22,2}$, $X_{23,2}$, and $X_{33,2}$ to bus 2. Bus 5 sends $X_{44,3}$ to bus 4. Bus 2 updates $\upsilon_{22,1}$, $\upsilon_{23,1}$ and $\upsilon_{33,1}$ by
\begin{align}
\upsilon_{ik,1}^{(t+1)} = \upsilon_{ik,1}^{(t)} - \alpha^{(t)}\left( X_{ik,2}-X_{ik,1}\right), i,k=2,3, i\leq k,
\label{ex2dualupdate1}
\end{align}
and bus 4 updates $\upsilon_{44,1}$ by
\begin{align}
\upsilon_{44,1}^{(t+1)} = \upsilon_{44,1}^{(t)} - \alpha^{(t)}\left( X_{44,3}-X_{44,2}\right).
\label{ex2dualupdate2}
\end{align}
The communication pattern is shown in Fig. \ref{fig:ring_dual}.

%---------------------------------------------------------------------------------------- 
\section{Simulation Results} \label{simulation}
\begin{table}[!t]
\renewcommand{\arraystretch}{1.3}
\caption{Normalized {CPU} time for distribution test feeders$^a$}
\label{tab:benchmark}
\centering
\begin{tabular}{c|c |c| c| c| c}
\hline\hline
\multirow{2}{1cm}{Number of buses} 	& \multirow{2}{1.3cm}{Centralized} 	& \multicolumn{2}{c|}{Cumulative} &  \multicolumn{2}{c}{Distributed}  \\ \cline{3-6}
				&						& Primal 	&	Dual		& Primal		& Dual \\
\hline
8 & 1.85	& 7.21	&	5.62	&	1.52	& 	1.00\\
34 & 298.68	& 37.79	&	33.89	&	1.94	& 	1.70\\
123 & --$^b$		& 143.39	&	126.48	&	2.24	& 	1.64\\
\hline\hline
\end{tabular}\\
\scriptsize
\begin{tabular}{p{7.5cm}}
$^a$ The CPU times are normalized by 0.0857s.\\
$^b$ The solver cannot be applied because of the out-of-memory problem.
\end{tabular}
\end{table}

\begin{table}[!t]
\renewcommand{\arraystretch}{1.3}
\caption{Normalized {CPU} time for IEEE Power transmission system benchmarks$^a$}
\label{tab:benchmark_trans}
\centering
\begin{tabular}{p{1cm}|p{1.3cm} |p{1.3cm}| p{1.3cm}|p{1cm}|p{1cm}}
\hline\hline
Number of buses 	& Centralized 	& Cumulative dual &  Distributed dual & iterations & Initial step size\\
\hline
14 & 5.38		& 5.38	&	1.00	& 1& 30\\
30 & 45.29	& 58.60	&	5.38	& 6 & 30\\
57 & 1696.79	& 49.08	&	4.28	& 4 & 30\\
118& --$^b$		& 704.46	&	13.51& 9& 300\\
\hline\hline
\end{tabular}
\scriptsize
\begin{tabular}{p{8.5cm}}
$^a$ The simulations for this problem set are done on MacBook Pro with 2.4 GHz Intel core i5 and 4 GB RAM. The CPU times are normalized by 0.1410s.\\
$^b$ The solver cannot be applied because of the out-of-memory problem.
\end{tabular}
\end{table}

\begin{figure}[!t]
	\begin{center}
		\includegraphics[width=0.53 \textwidth]{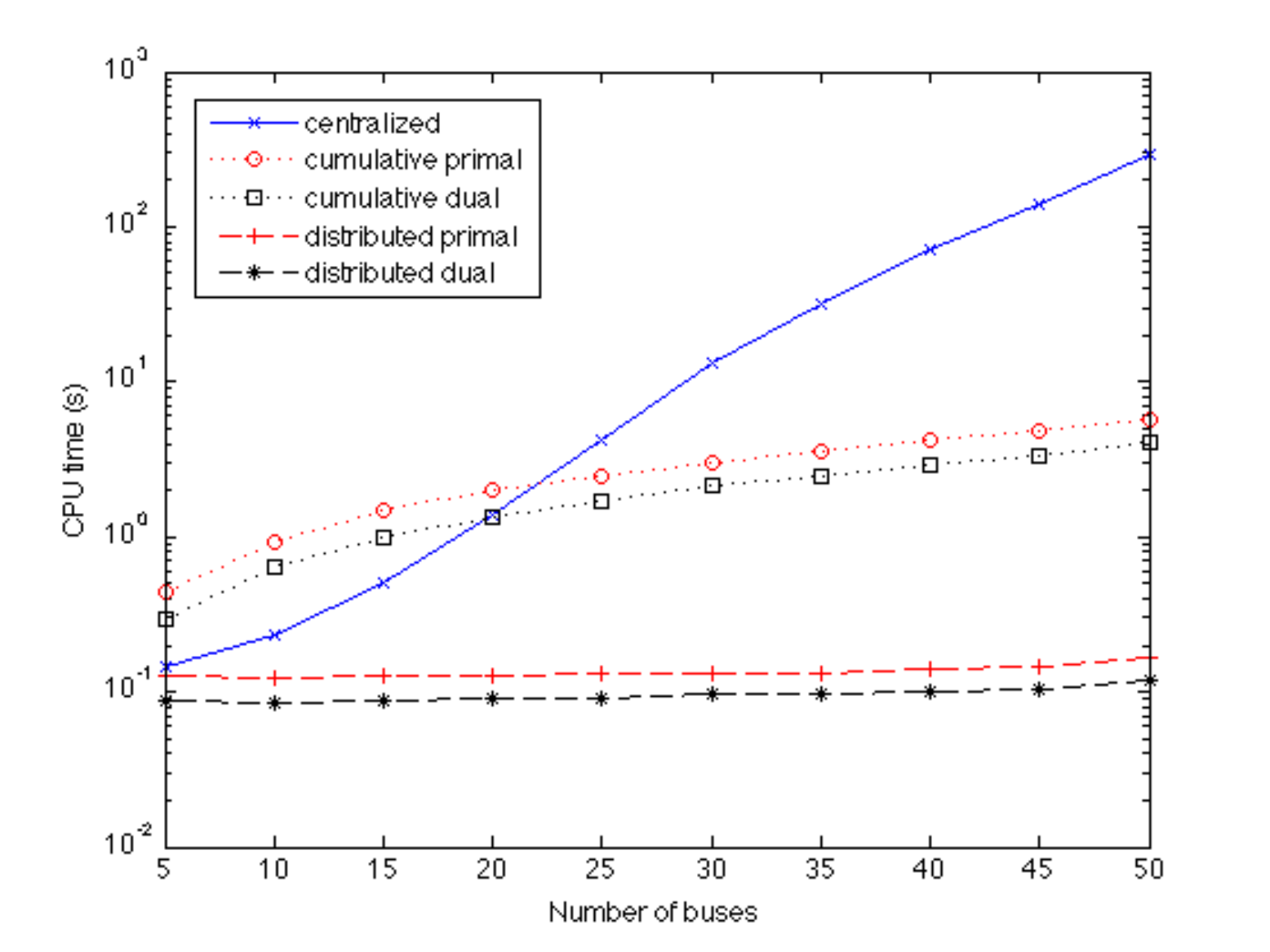}
	\end{center}
	\caption{CPU time of the various approaches on radial networks with bounded voltages}
  \label{plot_range}
\end{figure}

\begin{figure}[!t]
	\begin{center}
		\includegraphics[width=0.5 \textwidth]{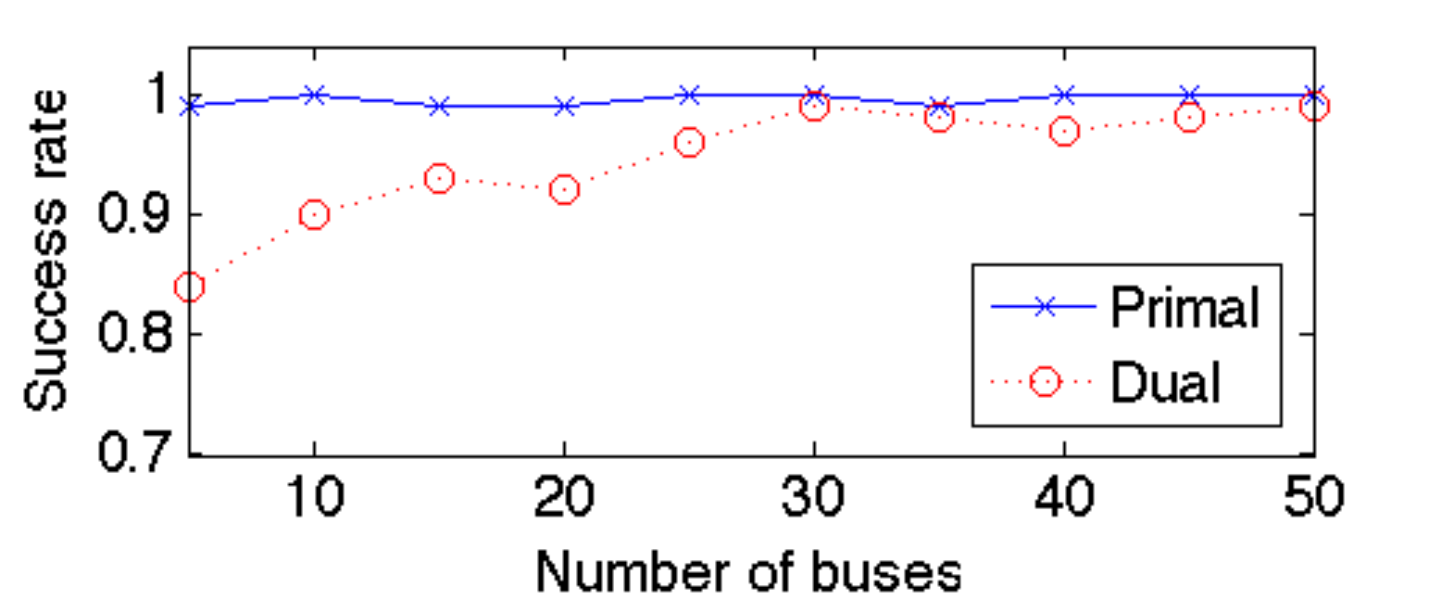}
	\end{center}
	\caption{Success rates of the primal and dual algorithms on radial networks with bounded voltages}
  \label{plot_rate}
\end{figure}

To evaluate the performance of the algorithms, we perform extensive simulations on various network settings. Since OPF is formulated as an SDP, the optimal solution can be computed in polynomial time by any popular SDP algorithms, e.g. IPM. Recall that the primal or dual algorithm aims to divide the original problem into smaller ones and to coordinate the subproblems, which of each can be solved by any SDP solver independently. The primal and dual algorithms perform coordination by indicating what problem data should be allocated to each subproblem and do simple calculations to update the shared terms (for the primal) and the prices (for the dual). When compared with those done by the SDP solver, the computation and ordination required solely by our algorithms are relatively far less stringent. As a whole, the bottleneck of computation should be at the SDP solver. In our simulation, we program the primal and dual algorithms in MATLAB and and solve each SDP with YALMIP \cite{yalmip} and SeDuMi \cite{sedumi}. To get rid of the dependence on the programming language and to simplify the comparison, we only count the CPU time spent on the SDP solver. Moreover, we can arrange the  subproblems to be solved in a single node or distribute them to different nodes in the network. For the former, we assume the problems are handled sequentially and we call it the cumulative approach. The latter, named as the distributed approach, addresses the subproblems in parallel. Without our algorithms, the (original) problem will be solved in its original form in a centralized manner. Here we compare the CPU times required for the SDP(s) among the centralized approach, (primal and dual) cumulative approaches, and (primal and dual) distributed approaches.

We run the simulations on Dell PowerEdge 2650 with 2 $\times$ 3.06GHz Xeon and 6GB RAM (except those for the transmission system benchmarks in Table \ref{tab:benchmark_trans}).\footnote{The results in Tables \ref{tab:benchmark} and \ref{tab:benchmark_trans} are normalized, and thus, they are comparable.}  In order to monitor the performance in each simulation run, we assemble the partial solutions (done by the subproblems) to form a complete one for the original problem and evaluate the corresponding objective value.\footnote{The assembly of partial solutions is not required in real implementation.} Our algorithms stop when the computed objective value falls in the range of $10^{-2}\times$ the global minimum. We assume that the dual algorithm is synchronized. In other words, all subproblems for the dual are solved in each iteration (but this is not required when implemented in real systems). The initial step size $\alpha^{(0)}$ is set to one and it is updated by $\alpha^{(t)}=\alpha^{(t-1)}/t, \forall t>0$. An algorithm is deemed successful if the stopping criterion is met in 100 iterations. 

We perform simulations on three problem sets; the first two focuses on tree-like networks while the last one is about transmission networks. The first problem set is  some distribution test feeder benchmarks \cite{feeder}. As the data set does not specify the cost function of power production/consumption, we create a problem instance by randomly generating the costs. To do this, we first select one node, e.g. node $i$, to be the power source node with $c_{i1}$ set randomly in the range $(0,10)$. For other node $k\neq i$, $c_{k1}$ set randomly in the range $(-10,0)$. We create 100 instances for each network. Table \ref{tab:benchmark} shows the averages of the normalized CPU times of the various approaches. All algorithms converge in 100 iterations for all instances.

The second problem set is  the $n$-bus radial network demonstrated in Section \ref{example}. For each instance, the root is the power source with a random cost selected in $(0,10)$ and each of the rest takes a random cost in  $(-10,0)$. For each node $i$, we specify a number $\xi$ in $(0.9,1.1)$ and set $\underline{V}_i=0.95\xi$ and $\overline{V}_i=1.05\xi$. For each line, the magnitudes of the conductance and susceptance are randomly assigned in $(0,10)$. We produce 100 instances for each $n$ and plot the average CPU times in logarithmic scale in Fig. \ref{plot_range}. 

From the simulation results for these two problem sets, both the primal and dual algorithms converge very fast for distribution networks. The CPU time for the centralized approach grows very fast with the size of the network.  The CPU time  grows roughly linearly for the cumulative approach while  it becomes independent of the size for the distributed approach. We define success rate as the fraction of the total number of simulation runs with stopping criterion met in 100 iterations, shown in Fig. \ref{plot_rate}. The success rate of the primal is almost $100\%$ for all tested network sizes. The dual fails to converge in 100 iterations for a small fraction of small networks but the success rate grows to almost $100\%$ with the network size. In general,the primal and dual algorithms are similar in performance but the dual requires a little bit less CPU time than the primal on the average.

The third problem set is some IEEE power transmission system benchmarks \cite{transmissionbenchmark}. As pointed out in \cite{zerogap}, these test cases have zero duality zap although they have network structures different from what we mention in Section \ref{sec:gap}. Table \ref{tab:benchmark_trans} shows the CPU times required, the iterations for convergence, and  the initial step sizes. The primal algorithm is not applied to this problem set and the reason will be given in the next section. For these transmission network topologies, the maximal cliques of the fill-in graphs are much irregular than those with tree-structured networks. There are many ways to construct the maximal cliques and different construction can result in different convergence speed. The study of the relationship between maximal clique construction and the algorithm performance is out of the scope of this paper. In this simulation, we randomly choose one maximal clique configuration and the step sizes are adjusted individually so as to have fast convergence. Nevertheless, the dual algorithm is more desirable than the centralized approach.

%----------------------------------------------------------------------------------------------
\section{Discussion} \label{discussion}

Both the primal and dual algorithms try to tackle the original problem by solving smaller subproblems but the ways to handle the information corresponding to the common partial solutions between subproblems are different. For the primal algorithm, if a common solution corresponds to a bus, that bus will compute its partial solution with the required information. For example, in  Fig.~\ref{fig:tree}, $W_{nn}$ for bus $n$ is common to all the subproblems. Bus $n$ computes $W_{nn}$ with its own $\underline{V}_n$, $\overline{V}_n$, and $M_{nn}$. Only the computed $W_{nn}$ is required to transmit to other buses which do not require bus $n$'s information. For the dual algorithm, each subproblem needs to acquire all its bus information, even for the common bus. Consider the example in Fig.~\ref{fig:tree} again, for $1\leq i \leq n-1$, node $i$ requires $\underline{V}_n$, $\overline{V}_n$, and $M_{nn}$ to solve its subproblem. Similarly, if the common solution is for a link, we can just elect any one of the connected bus to do the computation with the primal and dual algorithms. Therefore, the primal algorithm requires less information sharing between buses and it favors situations with sensitive bus information.

Our algorithms do not require an overlay of communication networks with different topology. From the examples in Section \ref{example}, all the communication is one-hop. A node needs to transmit data, e.g. $W,\lambda, X, \text{and } \upsilon$, to its neighbor nodes only. Communication links need to be built along with existing transmission lines only.

The primal algorithm is suitable for networks with tree structure while the dual can handle those with rings. In fact, the primal is not very efficient to update the partial solutions for (fill-in) edges, i.e. (\ref{edgere})--(\ref{fedgeim}), especially when their values are closed to the boundary of the feasible region. %When there is a common (fill-in) edge between subproblems, the connected nodes of the edge are also common to the subproblems too. 
When updating the variables, the bus variables, i.e. $W_{ii}$, are bounded but it is not the case for the edge variables, i.e. $W_{ik}, i\neq k$. When the step size is too large, (\ref{edgere})--(\ref{fedgeim}) may drive the current point out of the feasible region. If so, the obtained $\lambda$ cannot be used to form a subgradient, and thus (\ref{edgere})--(\ref{fedgeim}) fail. When reaching an infeasible point from a feasible one, we can always step back and choose a smaller step for an update again. When approaching the optimal solution which is located on the boundary for SDP, this step-back procedure is not very efficient. The dual algorithm does not have this problem when updating $\upsilon$ with (\ref{dualpriceupdate}). As mentioned, we can always constrain a feasible solution by averaging the shared variables.

Tree networks have a very nice property with respect to maximal clique construction; each branch with its attached nodes form a maximal clique. Hence, it is straightforward to decompose the problem into subproblems. However, when it comes to a more irregular network like those IEEE transmission system benchmarks, the number of ways to decompose the problem  grows with the size of the network. The performance of our algorithms also depends on how we form the maximal cliques and the initial step size needs to be adjusted accordingly. For problems with tree structure, the performance is easier to predict and the algorithms converge faster.

The dual algorithm is more resistant to communication delay than the primal. For the primal, an update of a shared variable requires $\lambda$ from all involved subproblems and thus synchronization is required. Delay of computing or transmitting an $\lambda$ from any subproblem can affect the whole algorithm proceed. On the other hand, an update of an $\upsilon$ requires the $\tilde{X}$ from two pre-associated subproblems according to the arrangement of the inequalities in \ref{slackequality}. Delay of computing or transmitting a particular $\tilde{X}$ can affect the update of some but not all the $\upsilon$. Thus the dual algorithm is asynchronous. Moreover, we can pair the variables in (\ref{slackequality}) into equalities differently and secretly whenever we start the dual algorithm. In some sense, the dual algorithm is more robust to attack stemmed from communication on the communication links.

%---------------------------------------------------------------------------------------- 
\section{Conclusion} \label{conclusion}
OPF is very important in planning the schedule of power generation. In the smart grid paradigm, more renewable energy sources will be incorporated into the system, especially in distribution networks, and the problem size will also grow tremendously. As problems with some special structures (e.g. trees for distribution networks) have a zero duality gap, we can find the optimal solution by solving the convex dual problem. In this paper, we propose the primal and  dual algorithms (with respect to the primal and dual decomposition techniques) to speed up the computation of the convexified  OPF problem. The problem is decomposed into smaller subproblems, each of which can be solved independently and effectively. The primal algorithm coordinates the subproblems by controlling the shared terms (related to electricity resources) while the dual one manages them by updating the prices. From the simulation results for tree-structure problems,  the computation time grows linearly with the problem size if we solve the decomposed problem in a central node with our algorithms. The computation time becomes independent of the problem size when the subproblems are solved in parallel in different nodes. Even without nice network structure such as a tree, the dual algorithm outperforms the centralized approach without decomposition. Therefore, the primal and dual algorithms are excellent in addressing OPF, especially for distribution networks. In future, we will improve the algorithm by incorporate more constraints into the OPF problem and move to nonlinear objective functions.

% or


%\appendices
%% you can choose not to have a title for an appendix
%% if you want by leaving the argument blank
%\section{}

%\section*{Acknowledgment}


\ifCLASSOPTIONcaptionsoff
  \newpage
\fi



\bibliographystyle{IEEEtran}
\bibliography{IEEEabrv}

\begin{thebibliography}{10}
\providecommand{\url}[1]{#1}
\csname url@samestyle\endcsname
\providecommand{\newblock}{\relax}
\providecommand{\bibinfo}[2]{#2}
\providecommand{\BIBentrySTDinterwordspacing}{\spaceskip=0pt\relax}
\providecommand{\BIBentryALTinterwordstretchfactor}{4}
\providecommand{\BIBentryALTinterwordspacing}{\spaceskip=\fontdimen2\font plus
\BIBentryALTinterwordstretchfactor\fontdimen3\font minus
  \fontdimen4\font\relax}
\providecommand{\BIBforeignlanguage}[2]{{%
\expandafter\ifx\csname l@#1\endcsname\relax
\typeout{** WARNING: IEEEtran.bst: No hyphenation pattern has been}%
\typeout{** loaded for the language `#1'. Using the pattern for}%
\typeout{** the default language instead.}%
\else
\language=\csname l@#1\endcsname
\fi
#2}}
\providecommand{\BIBdecl}{\relax}
\BIBdecl

\bibitem{smartgrid}
P.~P. Varaiya, F.~F. Wu, and J.~W. Bialek, ``Smart operation of smart grid:
  Risk-limiting dispatch,'' \emph{Proc. {IEEE}}, vol.~99, pp. 40--57, 2011.

\bibitem{dcopf}
B.~Stott, J.~Jardim, and O.~Alsac, ``Dc power flow revisited,'' \emph{{IEEE}
  Trans. Power Syst.}, vol.~24, pp. 1290--1300, 2009.

\bibitem{zerogap}
J.~Lavaei and S.~H. Low, ``Zero duality gap in optimal power flow problem,''
  \emph{{IEEE} Trans. Power Syst.}, in press.

\bibitem{geometry}
B.~Zhang and D.~Tse, ``Geometry of feasible injection region of power
  networks,'' \emph{To appear in proc. Allerton}, 2011.

\bibitem{network_zero}
S.~Sojoudi and J.~Lavaei, ``Network topologies guaranteeing zero duality gap
  for optimal power flow problem,'' \emph{Submitted to {IEEE} Trans. Power
  Syst.}, 2011.

\bibitem{coneprogram}
R.~A. Jabr, ``Radial distribution load flow using conic programming,''
  \emph{{IEEE} Trans. Power Syst.}, vol.~21, pp. 1458--1459, 2006.

\bibitem{kelly}
F.~Kelly, A.~Maulloo, and D.~Tan, ``Rate control in communication networks:
  shadow prices, proportional fairness and stability,'' \emph{Journal of the
  Operational Research Society}, vol.~49, pp. 237--252, 1998.

\bibitem{nonlinear}
D.~P. Bertsekas, \emph{Nonlinear programming}, 2nd~ed.\hskip 1em plus 0.5em
  minus 0.4em\relax Belmont, MA, USA: Athena Scientific, 1999.

\bibitem{chiang}
M.~Chiang, S.~H. Low, A.~R. Calderbank, and J.~C. Doyle, ``Layering as
  optimization decomposition: A mathematical theory of network architectures,''
  \emph{Proc. {IEEE}}, vol.~95, pp. 255--312, Jan. 2007.

\bibitem{pdipm}
S.~Wright, \emph{Primal-dual interior-point methods}.\hskip 1em plus 0.5em
  minus 0.4em\relax Philadelphia, PA, USA: Society for Industrial and Applied
  Mathematics, 1997.

\bibitem{MCS}
R.~E. Tarjan and M.~Yannakakis, \emph{Simple linear-time algorithms to test
  chordality of graphs, test acyclicity of hypergraphs, and selectively reduce
  acyclic hypergraphs}.\hskip 1em plus 0.5em minus 0.4em\relax Philadelphia,
  PA, USA: Society for Industrial and Applied Mathematics, July 1984, vol.~13.

\bibitem{elimination}
D.~R. Fulkerson and O.~A. Gross, ``Incidence matrices and interval graph,''
  \emph{Pacific J. Math.}, vol.~15, no.~3, pp. 835--855, 1965.

\bibitem{FIC}
R.~E. Neapolitan, \emph{Probabilistic reasoning in expert systems: theory and
  algorithms}.\hskip 1em plus 0.5em minus 0.4em\relax New York, NY, USA: John
  Wiley \& Sons, Inc., 1990.

\bibitem{maximalclique}
E.~A. Akkoyunlu, ``The enumeration of maximal cliques of large graphs,''
  \emph{SIAM Journal on Computing}, vol.~2, pp. 1--6, 1973.

\bibitem{parallelSDP}
K.~Nakata, M.~Yamashita, K.~Fujisawa, and M.~Kojima, ``A parallel primal-dual
  interior-point method for semidefinite programs using positive definite
  matrix completion,'' \emph{Parallel Comput.}, vol.~32, pp. 24--43, Jan. 2006.

\bibitem{sdp1}
M.~Fukuda, M.~Kojima, K.~Murota, and K.~Nakata, ``Exploiting sparsity in
  semidefinite programming via matrix completion {I}: General framework,''
  \emph{SIAM J. on Optimization}, vol.~11, pp. 647--674, March 2000.

\bibitem{sdp2}
K.~Nakata, K.~Fujisawa, M.~Fukuda, M.~Kojima, and K.~Murota, ``Exploiting
  sparsity in semidefinite programming via matrix completion {II}:
  implementation and numerical results,'' \emph{Mathematical Programming},
  vol.~95, pp. 303--327, 2003.

\bibitem{partial}
R.~Grone, C.~R. Johnson, E.~M. Sa, and H.~Wolkowicz, ``Positive definite
  completions of partial hermitian matrices,'' \emph{Linear Algebra and Its
  Applications}, vol.~58, pp. 109--124, 1984.

\bibitem{convex}
S.~Boyd and L.~Vandenberghe, \emph{Convex Optimization}.\hskip 1em plus 0.5em
  minus 0.4em\relax New York, NY, USA: Cambridge University Press, 2004.

\bibitem{yalmip}
J.~L{\"o}fberg, ``{YALMIP : A toolbox for modeling and optimization in
  MATLAB},'' in \emph{Proc. International Symposium on Computer Aided Control
  Systems Design}, Sep. 2004, pp. 284--289.

\bibitem{sedumi}
J.~F. Sturm, ``Using sedumi 1.02, a matlab toolbox for optimization over
  symmetric cones,'' 1998.

\bibitem{feeder}
W.~H. Kersting, ``Radial distribution test feeders,'' in \emph{Proc. IEEE Power
  Engineering Society Winter Meeting}, vol.~2, 2001, pp. 908--912.

\bibitem{transmissionbenchmark}
\BIBentryALTinterwordspacing
University of washington, power systems test case archive. [Online]. Available:
  \url{http://www.ee.washington.edu/research/pstca}
\BIBentrySTDinterwordspacing

\end{thebibliography}
%\begin{thebibliography}{1}
%\end{thebibliography}

\end{document}